\documentclass[12pt]{amsart}

\usepackage[margin=1in]{geometry}
\usepackage{amsmath,amsthm,amssymb}
\usepackage{dsfont}
\usepackage{amsmath,amscd}%CD package
\usepackage{mathrsfs}%mathscr
\usepackage{color}%hyperlink
\usepackage{hyperref}%hyperlink
\usepackage{enumitem}%enumerate label
\usepackage{array}%tabular width
\usepackage{graphicx}
\usepackage{xcolor}

\usepackage{soul}%strokeout

\hypersetup{
    colorlinks=true, % make the links colored
    linkcolor=blue, % color TOC links in blue
    urlcolor=red, % color URLs in red
    linktoc=all % 'all' will create links for everything in the TOC
}

\input xypic 
\input xy
\xyoption{all}

\newcommand{\Q}{\mathbb{Q}}
\newcommand{\R}{\mathbb{R}}
\newcommand{\C}{\mathbb{C}}

\newcommand{\PP}{\mathbb{P}}
\newcommand{\Z}{\mathbb{Z}}

\newcommand{\G}{\mathcal{G}}
\newcommand{\map}{\textrm{Map}}
\newcommand{\qqed}{\hfill\Box}

\theoremstyle{plain}
\newtheorem{thm}{Theorem}[section]
\newtheorem{prop}[thm]{Proposition}
\newtheorem{lemma}[thm]{Lemma}
\newtheorem{lemma'}{Lemma}
\newtheorem{cor}[thm]{Corollary}

\theoremstyle{definition}

\newcommand{\namedright}[3]{\ensuremath{#1\stackrel{#2}
 {\longrightarrow}#3}}
\newcommand{\nameddright}[5]{\ensuremath{#1\stackrel{#2}
 {\longrightarrow}#3\stackrel{#4}{\longrightarrow}#5}}
\newcommand{\namedddright}[7]{\ensuremath{#1\stackrel{#2}
 {\longrightarrow}#3\stackrel{#4}{\longrightarrow}#5
  \stackrel{#6}{\longrightarrow}#7}}

\newcommand{\larrow}{\relbar\!\!\relbar\!\!\rightarrow}
\newcommand{\llarrow}{\relbar\!\!\relbar\!\!\larrow}

\newcommand{\lnameddright}[5]{\ensuremath{#1\stackrel{#2}
 {\larrow}#3\stackrel{#4}{\larrow}#5}}

\newcommand{\llnamedright}[3]{\ensuremath{#1\stackrel{#2}
 {\llarrow}#3}}
\newcommand{\llnameddright}[5]{\ensuremath{#1\stackrel{#2}
 {\llarrow}#3\stackrel{#4}{\llarrow}#5}}

\newcolumntype{C}[1]{>{\centering\let\newline\\\arraybackslash\hspace{0pt}}m{#1}}

\begin{document}

\title[Suspension of 4-manifolds]{The suspension of a $4$-manifold and its applications} 
\author{Tseleung So}
\address{Department of Mathematics, Univeristy of Western Ontario, London ON, N6A 5B7, Canada}
\email{tso28@uwo.ca}
\author{Stephen Theriault}
\address{Mathematical Sciences, University
         of Southampton, Southampton SO17 1BJ, United Kingdom}
\email{S.D.Theriault@soton.ac.uk}

\subjclass[2010]{Primary 57N65, 55P40, Secondary 81T13.} 
\keywords{non-simply-connected four-manifold, suspension, homotopy type, gauge group}

\maketitle

\begin{abstract}
Let $M$ be a smooth, orientable, closed, connected $4$-manifold 
and suppose that~$H_1(M;\mathbb{Z})$ is finitely generated and has no $2$-torsion. We give a 
homotopy decomposition of the suspension of $M$ in terms of spheres, Moore spaces 
and $\Sigma\mathbb{C}P^{2}$. This is used to calculate any reduced 
generalized cohomology theory of $M$ as a 
group and to determine the homotopy types of certain current groups and gauge groups.
\end{abstract}

\section{Introduction} 

Let $M$ be a smooth, orientable, closed, connected $4$-manifold. This implies by Morse 
theory that~$M$ has a $CW$-structure with one $4$-cell. Suppose that $H_{1}(M;\mathbb{Z})$ 
is finitely generated and has no $2$-torsion. Specifically, assume that: 
\begin{equation} 
  \label{Mhyp} 
  \begin{split}
       \bullet & \ H_{1}(M;\mathbb{Z})\cong\mathbb{Z}^{m}\oplus 
                   \bigoplus_{j=1}^{n}\mathbb{Z}/b_{j}\mathbb{Z}; \\   
       \bullet & \ \mbox{each $b_{j}$ is a prime power, where the prime is odd}.
  \end{split} 
\end{equation} 
From~(\ref{Mhyp}), by Poincar\'{e} Duality, the integral homology of $M$ is: 
\begin{equation} 
\label{MPD} 
\begin{array}{c|c}
i	&H_i(M;\mathbb{Z})\\
\hline
0	&\Z\\[1mm]
1	&\Z^m\oplus\bigoplus^n_{j=1}\Z/b_j\Z\\[2mm]
2	&\Z^d\oplus\bigoplus^n_{j=1}\Z/b_j\Z\\[2mm]
3	&\hspace{3mm}\Z^m\\[1mm]
4	&\Z\\[1mm]
\geq5	&0
\end{array}
\end{equation} 
where $d\geq 0$ can be any integer. Our main theorem identifies the homotopy 
type of $\Sigma M$. 

\begin{thm} 
   \label{suspM} 
   Let $M$ be a smooth, orientable, closed, connected $4$-manifold 
   and suppose that~$H_1(M;\mathbb{Z})$ is finitely generated and has no $2$-torsion. If $M$ is Spin then there is a homotopy equivalence 
   \[\Sigma M\simeq\bigg(\bigvee_{i=1}^{m} (S^{2}\vee S^{4})\bigg)\vee 
          \bigg(\bigvee_{j=1}^{n} (P^{3}(b_{j})\vee P^{4}(b_{j}))\bigg)\vee 
          \bigg(\bigvee_{k=1}^{d} S^{3}\bigg)\vee S^{5}.\] 
   If $M$ is non-Spin then there is a homotopy equivalence 
   \[\Sigma M\simeq\bigg(\bigvee_{i=1}^{m} (S^{2}\vee S^{4})\bigg)\vee 
          \bigg(\bigvee_{j=1}^{n} (P^{3}(b_{j})\vee P^{4}(b_{j}))\bigg)\vee  
          \bigg(\bigvee_{k=1}^{d-1} S^{3}\bigg)\vee\Sigma\mathbb{C}P^{2}.\] 
\end{thm} 

In fact, Theorem~\ref{suspM} is a special case of a more general result about the suspension 
of 4-dimensional $CW$-complexes whose cohomology satisfies Poincar\'{e} Duality 
and has no $2$-torsion (see Theorem~\ref{suspM CW case}). Such a classification 
fits into a long history of classifying $CW$-compexes with cells occurring in a small number of 
consecutive dimensions, with contributions, for example, by Whitehead~\cite{W1,W2}, Chang~\cite{C}, 
Baues and Hennes~\cite{BH}, Baues and Drozd~\cite{BD} and Pan and Zhu~\cite{PZ}.
Apart from~\cite{W2}, these classifications occur in the stable range; the classification in 
Theorem~\ref{suspM CW case} notably occurs unstably. 

A key aspect of Theorem~\ref{suspM} is that the suspension of $M$ involves 
only three types of spaces: spheres, Moore spaces and $\Sigma\mathbb{C}P^{2}$. Each 
is simple and characterizes a cohomological property: a sphere corresponds to an isolated 
$\mathbb{Z}$ summand, a Moore space corresponds to a torsion summand, and a 
$\Sigma\mathbb{C}P^{2}$ corresponds to two $\mathbb{Z}$ summands connected by 
the Steenrod operation $Sq^{2}$. The hypothesis that only odd torsion in cohomology 
is allowed is necessary to achieve this. For example, the suspension of $S^1\times\R\PP^3$ 
is homotopy equivalent to $S^2\vee\Sigma\R\PP^3\vee\Sigma^2\R\PP^3$ which does not split 
as in Theorem~\ref{suspM} since $\Sigma\R\PP^3$ is indecomposable. The list of 
indecomposable wedge summands at the 
prime $2$ would therefore be much more complex. 

The simple description of $\Sigma M$ in Theorem~\ref{suspM} is advantageous. It
implies that the homotopy type of $\Sigma M$ is completely determined by only two 
properties: (i) whether $M$ is Spin or not and~(ii) $H_{\ast}(M;\mathbb{Z})$ 
(or equivalently,~$H^{\ast}(M;\mathbb{Z})$). 

Interestingly, while suspending a manifold loses all the geometry, it does give access to 
many other properties. Theorem~\ref{suspM} is applied in three different contexts: to determine 
any reduced generalized cohomology theory of $M$, to determine the homotopy 
type of certain current groups associated to $M$, and to determine the homotopy type 
of certain gauge groups associated to~$M$. These applications are discussed in detail 
in Section~\ref{sec:apps}. 

To prove Theorem~\ref{suspM} new methods are developed that use homology and cohomology to detect 
whether certain maps are null homotopic. This generalizes Neisendorfer's work in defining 
and determining the mod-$p^{r}$ Hopf invariant~\cite{neisendorfer10}.

\section{Preliminary information on Moore spaces} 
\label{sec:moore} 

This section records some information on the homotopy groups of Moore 
spaces which will be needed later. For $m\geq 2$ and $k\geq 2$, 
the \emph{mod-$k$ Moore space} $P^{m}(k)$ of dimension $m$ is the homotopy 
cofibre of the degree $k$ map on $S^{m-1}$. Notice that $\Sigma P^{m}(k)\simeq P^{m+1}(k)$. 

\begin{lemma} 
   \label{pi3Moore} 
   If $p$ is an odd prime and $r\geq 1$ then $\pi_{3}(P^{3}(p^{r}))\cong\mathbb{Z}/p^{r}\mathbb{Z}$. 
\end{lemma} 

\begin{proof} 
Consider the homotopy fibration 
\(\nameddright{F^{3}(p^{r})}{}{P^{3}(p^{r})}{q}{S^{3}}\) 
where $q$ is the pinch map to the top cell. This induces an exact sequence 
\[\namedddright{[S^{3},\Omega S^{3}]}{}{[S^{3},F^{3}(p^{r})]}{}{[S^{3},P^{3}(p^{r})]} 
       {q_{\ast}}{[S^{3},S^{3}]}.\] 
At odd primes, $\pi_3(\Omega S^3)\cong0$. Since $P^{3}(p^{r})$ is rationally trivial 
and $\pi_3(S^3)\to\pi_3(S^3)\otimes\Q$ is injective, 
any composite  
\(\nameddright{S^{3}}{f}{P^{3}(p^{r})}{q}{S^{3}}\) 
must have degree zero. Hence $q_{\ast}=0$. Thus, by exactness, 
$\pi_{3}(F^{3}(p^{r}))\cong\pi_{3}(P^{3}(p^{r}))$. 

To complete the proof it is now equivalent to show that 
$\pi_{2}(\Omega F^{3}(p^{r}))\cong\mathbb{Z}/p^{r}\mathbb{Z}$. 
For $m\geq 1$, let $S^{2m+1}\{p^{r}\}$ be the homotopy fibre of the degree $p^{r}$ 
map on $S^{2m+1}$. In particular, $S^{2m+1}\{p^{r}\}$ is $(2m-1)$-connected. 
By~\cite[Proposition 14.2]{neisendorfermem} there is a homotopy equivalence 
\[\Omega F^{3}(p^{r})\simeq S^{1}\times\big(\prod_{j=1}^{\infty} S^{2p^{j}-1}\{p^{r+1}\}\big)\times 
       \Omega R^{3}(p^{r})\] 
where $R^{3}(p^{r})$ is a wedge of mod-$p^{r}$ Moore spaces consisting of a single 
copy of $P^{4}(p^{r})$ and all other wedge summands being at least $3$-connected. 
In particular, for $R^{3}(p^{r})$, by the Hilton-Milnor Theorem there is an isomorphism 
$\pi_{3}(R^{3}(p^{r}))\cong\pi_{3}(P^{4}(p^{r}))$. Further, the Hurewicz homomorphism  
implies that $\pi_{3}(P^{4}(p^{r}))\cong H_{3}(P^{4}(p^{r}))\cong\mathbb{Z}/p^{r}\mathbb{Z}$. 
Returning to the decomposition of $\Omega F^{3}(p^{r})$, since each space $S^{2p^{j}-1}\{p^{r+1}\}$ 
is at least $3$-connected, we obtain~$\pi_{2}(\Omega F^{3}(p^{r}))\cong\pi_{2}(\Omega R^{3}(p^{r}))$ and we have just seen 
that $\pi_{2}(\Omega R^{3}(p^{r}))\cong\mathbb{Z}/p^{r}\mathbb{Z}$. 
\end{proof} 

\begin{lemma}\cite[Lemma 3.3]{so16}
\label{lemma_pi4 P3, pi4 P4 are trivial}
If $p$ is an odd prime and $r\geq 1$ then $\pi_4(P^3(p^r))\cong 0$ and~$\pi_4(P^4(p^r))\cong 0$.~$\qqed$ 
\end{lemma} 

\begin{lemma} 
   \cite[Corollary 6.6]{neisendorfermem} 
   \label{Mooreskel} 
   Let $p$ be an odd prime, $s,t\geq 1$ and $m,n\geq 2$. Then there is a homotopy equivalence 
   \[P^{m}(p^{s})\wedge P^{n}(p^{t})\simeq P^{m+n-1}(p^{min(s,t)})\vee P^{m+n}(p^{min(s,t)}).\]
\end{lemma} 
\vspace{-0.8cm}$\qqed$\bigskip 

%\begin{proof}
%Assume $s\leq t$. Applying the functor $P^m(p^s)\wedge-$ to the cofibration $S^{n-1}\overset{p^t}{\to}S^{n-1}\to %P^n(p^t)$ gives a cofibration
%\[
%P^{m+n-1}(p^s)\overset{p^t}{\to}P^{m+n-1}(p^s)\to P^{m}(p^s)\wedge P^n(p^t)\to P^{m+n}(p^s)\overset{p^s}{\to}%P^{m+n}(p^s).
%\]
%The degree maps $p^t$ are zero maps so $P^m(p^s)\wedge P^n(p^t)\simeq P^{m+n-1}(p^s)\vee P^{m+n}(p^s)$.
%\end{proof}
%}

\begin{lemma} 
   \label{pi3mooresmash} 
   Let $p$ be an odd prime and $s,t\geq 1$. Then 
   $\pi_{3}(\Sigma P^{2}(p^{s})\wedge P^{2}(p^{t}))\cong\mathbb{Z}/p^{min(s,t)}\mathbb{Z}$. 
\end{lemma} 

\begin{proof} 
By Lemma~\ref{Mooreskel} and for dimensional reasons there are isomorphisms 
\[\pi_{3}(\Sigma P^{2}(p^{s})\wedge P^{2}(p^{t}))\cong 
     \pi_{3}(P^{4}(p^{min(s,t)})\vee P^{5}(p^{min(s,t)})\cong\pi_{3}(P^{4}(p^{min(s,t)})).\]
Since $P^{4}(p^{min(s,t)})$ is $2$-connected, by the Hurewicz Theorem there are 
isomorphisms 
\[\pi_{3}(P^{4}(p^{\min(s,t)}))\cong H_{3}(P^{4}(p^{min(s,t)});\mathbb{Z})\cong 
   \mathbb{Z}/p^{min(s,t)}\mathbb{Z}.\] 
\end{proof}

\section{A homological test for a null homotopy I} 
\label{sec:test1} 

In the next two sections we give homological and cohomological criteria determining 
when certain maps are null homotopic. These maps are from $S^{3}$ or $P^{3}(p^{r})$ 
into a wedge~$\bigvee_{i=1}^{m} P^{3}(p^{r_{i}})$. So the material in this section and the 
next focus on $3$-dimensional Moore spaces. 

In what follows we will use the terms ``homotopy fibration diagram" and 
``homotopy cofibration diagram". To explain these, recall that there is a standard construction 
that turns any continuous, pointed map 
\(f\colon\namedright{X}{}{Y}\) 
that is a surjection on path-components into a fibration, in the sense that 
$f$ factors as $p\circ\phi$ where 
\(\phi\colon\namedright{X}{}{X'}\) 
is a homotopy equivalence and 
\(p\colon\namedright{X'}{}{Y}\) 
is a fibration  (see, for example, \cite[Theorem 7.1.14]{Se}). The homotopy fibre of $f$ is 
the fibre of $p$. As in~\cite[Section 7.6]{Se}, a homotopy commutative square 
\begin{equation} 
\label{Selickstart} 
\diagram 
      W\rto^-{g'}\dto^{f'} & X\dto^{f} \\ 
      Y\rto^-{g} & Z 
  \enddiagram 
\end{equation} 
is equivalent up to homotopy to a strictly commutative square in which the horizontal maps 
are fibrations. This induces a map between fibres, that is, a map between the homotopy fibres 
of $g'$ and $g$. It is notable that while the homotopy types of the fibres are determined by the 
homotopy classes of $g'$ and $g$, the homotopy class of the induced map is not determined 
by the homotopy classes of $f$ and $f'$. However, the induced map $\gamma$ can be chosen 
via the standard construction above so that there is a homotopy commutative diagram of fibration sequences 
\[\diagram 
      \Omega X\rto^-{\partial'}\dto^{\Omega f} & F'\rto\dto^{\gamma} & W\rto^-{g}\dto^{f'} & X\dto^{f} \\ 
      \Omega Y\rto^-{\partial} & F\rto & Y\rto^-{g}  & Z. 
  \enddiagram\] 
Further, this diagram could be extended vertically as well, as in~\cite[Thoerem 7.6.2]{Se}, to produce 
a homotopy commutative diagram in which each consecutive pair of horizontal maps and each 
consecutive pair of vertical maps is a homotopy fibration. Any such diagram originating from 
the square~(\ref{Selickstart}) and extending via homotopy fibrations horizontally or vertically 
in this manner is called a \emph{homotopy fibration diagram}.  A \emph{homotopy cofibration diagram} 
is defined dually.  

In general, let 
\(i_{1}\colon\namedright{\Sigma X}{}{\Sigma X\vee\Sigma Y}\)  
and 
\(i_{2}\colon\namedright{\Sigma Y}{}{\Sigma X\vee\Sigma Y}\) 
be the inclusions of the left and right wedge summands respectively. Let 
\[[i_{1},i_{2}]\colon\namedright{\Sigma X\wedge Y}{}{\Sigma X\vee\Sigma Y}\] 
be the Whitehead product of $i_{1}$ and $i_{2}$. 

Let $r,s,t$ be positive integers such that $s,t\geq r$. Then 
\[H^2(P^3(p^s);\mathbb{Z}/p^{r}\mathbb{Z})\cong H^2(P^3(p^t);\mathbb{Z}/p^{r}\mathbb{Z})\cong 
     \mathbb{Z}/p^{r}\mathbb{Z}.\] 
Let $u_s$ and $u_t$ be the generators of $H^2(P^3(p^s);\mathbb{Z}/p^{r}\mathbb{Z})$ and $H^2(P^3(p^t);\mathbb{Z}/p^{r}\mathbb{Z})$ respectively. Then 
$H^2(P^3(p^s)\times P^3(p^t);\mathbb{Z}/p^{r}\mathbb{Z})$ is generated by 
$u_{s}\otimes 1$ and $1\otimes u_{t}$.

\begin{lemma} 
   \label{H4lemma} 
   Let $p$ be a prime and let $s$ and $t$ be integers such that $s,t\geq 1$. 
   Then there is an isomorphism 
   \[H^{4}(P^{2}(p^{s})\times P^{2}(p^{t});\mathbb{Z}/p^{\min(s,t)}\mathbb{Z})\cong\mathbb{Z}/p^{\min(s,t)}\mathbb{Z}\]
   and $u_{s}\cup u_{t}$ is a generator.
\end{lemma} 

\begin{proof}  
One case of the K\"{u}nneth Theorem (see, for example, \cite[Theorem 3.15]{hatcher}) is as follows. 
If $X$ and $Y$ are $CW$-complexes, $R$ is a ring, and $H^{k}(Y;R)$ is a finitely generated 
$R$-module for all $k$ then the cross product 
\(\namedright{H^{\ast}(X;R)\otimes_{R} H^{\ast}(Y;R)}{}{H^{\ast}(X\times Y;R)}\) 
is a ring isomorphism. In our case, if $r=\min(s,t)$ then both 
$H^{\ast}(P^{2}(p^{s});\mathbb{Z}/p^{r}\mathbb{Z})$ 
and $H^{\ast}(P^{2})(p^{t});\mathbb{Z}/p^{r}\mathbb{Z})$ are finitely generated free 
$\mathbb{Z}/p^{r}\mathbb{Z}$-modules. Therefore, by the K\"{u}nneth Theorem, there are isomorphisms 
\[\begin{split} 
   H^{4}(P^{2}(p^{s})\times P^{2}(p^{t});\mathbb{Z}/p^{r}\mathbb{Z}) &  
        \cong H^{2}(P^{2}(p^{s});\mathbb{Z}/p^{r}\mathbb{Z})\otimes 
                  H^{2}(P^{2}(p^{t});\mathbb{Z}/p^{r}\mathbb{Z}) \\  
        & \cong \mathbb{Z}/p^{r}\mathbb{Z}\otimes\mathbb{Z}/p^{r}\mathbb{Z}\cong 
                  \mathbb{Z}/p^{r}\mathbb{Z} 
  \end{split}\] 
and $u_{s}\cup u_{t}$ is a generator.
\end{proof} 

Propositions~\ref{S3cuptriv} and~\ref{Moorecuptriv} give useful tests 
for when a certain map is null homotopic. 

\begin{prop} 
\label{S3cuptriv}
Let $p$ be an odd prime and $s, t\geq 1$. Let 
$f:S^{3}\to\Sigma P^2(p^s)\wedge P^2(p^t)$ 
be a map and let $C$ be the homotopy cofiber of the composite
\[
S^{3}\overset{f}{\longrightarrow}\Sigma P^2(p^s)\wedge P^2(p^t) 
    \overset{[\imath_1,\imath_2]}{\longrightarrow}P^3(p^s)\vee P^3(p^t).
\] 
The following are equivalent: 
\begin{itemize} 
   \item[(a)] the map $f$ is null homotopic; 
   \item[(b)] \(\namedright{H^{3}(\Sigma P^{2}(p^{s})\wedge P^{2}(p^{t});\mathbb{Z}/p^{min(s,t)}\mathbb{Z})} 
                        {f^{\ast}}{H^{3}(S^3;\mathbb{Z}/p^{min(s,t)}\mathbb{Z})}\) 
                  is the zero map; 
   \item[(c)] all cup products in $\widetilde{H}^*(C;\Z/p^{min(s,t)}\Z)$ are zero. 
\end{itemize} 
\end{prop} 

\begin{proof} 
(a) $\Leftrightarrow$ (b). Let $u=min(s,t)$ and consider the following string of isomorphisms: 
\[\begin{split} 
      \pi_{3}(\Sigma P^{2}(p^{s})\wedge P^{2}(p^{t})) 
           & \cong H_{3}(\Sigma P^{2}(p^{s})\wedge P^{2}(p^{t});\mathbb{Z}) \\ 
           & \cong H_{3}(P^{4}(p^{u})\vee P^{5}(p^{u});\mathbb{Z}) \\ 
           & \cong H_{3}(P^{4}(p^{u})\vee P^{5}(p^{u});\mathbb{Z}/p^{u}\mathbb{Z}) \\ 
           & \cong H^{3}(P^{4}(p^{u})\vee P^{5}(p^{u});\mathbb{Z}/p^{u}\mathbb{Z}) \\ 
           & \cong H^{3}(\Sigma P^{2}(p^{s})\wedge P^{2}(p^{r});\mathbb{Z}/p^{u}\mathbb{Z}) 
   \end{split}\]
The first isomorphism is due to the Hurewicz Theorem because 
$\Sigma P^{2}(p^{s})\wedge P^{2}(p^{t})$ is $2$-connected. The second isomorphism 
holds by Lemma~\ref{Mooreskel}. The third isomorphism holds since~$H_{3}(P^{4}(p^{u})\vee P^{5}(p^{u});\mathbb{Z})\cong 
     H_{3}(P^{4}(p^{u});\mathbb{Z})\cong\mathbb{Z}/p^{u}\mathbb{Z}$
and changing homology coefficients from $\mathbb{Z}$ to $\mathbb{Z}/p^{u}\mathbb{Z}$ 
induces an isomorphism here. The fourth isomorphism holds by the Universal 
Coefficient Theorem. The fifth isomorphism holds by Lemma~\ref{Mooreskel}. Observe 
that under these isomorphisms the map 
\(\namedright{S^{3}}{f}{\Sigma P^{2}(p^{s})\wedge P^{2}(p^{t})}\) 
is sent to
\[\namedright{H^{3}(\Sigma P^{2}(p^{s})\wedge P^{2}(p^{t});\mathbb{Z}/p^{u}\mathbb{Z})} 
      {f^{\ast}}{H^{3}(S^{3};\mathbb{Z}/p^{u}\mathbb{Z})}.\]
Thus~$f$ is null homotopic if and only if $f^{\ast}=0$ in degree~$3$ mod-$p^{u}$ 
cohomology.~\medskip 

\noindent 
(a) $\Rightarrow$ (c). If $f$ is null homotopic then $C\simeq P^{3}(p^{s})\vee P^{3}(p^{t})\vee S^{4}$ 
is a suspension, so all cup products in $\widetilde{H}^*(C;\Z/p^{u}\Z)$ are zero. 
\medskip 

\noindent 
(c) $\Rightarrow$ (b). 
Consider the homotopy cofibration diagram 
\begin{equation}\label{CD_compare C_f and P^3xP^3}
\xymatrix{
S^3\ar@{=}[d]\ar[r]^-{f}						&\Sigma P^2(p^s)\wedge P^2(p^t)\ar[d]^-{[\imath_1,\imath_2]}\ar[r]	&C_f\ar[d]\\
S^3\ar[d]\ar[r]^-{[\imath_1,\imath_2]\circ f}	&P^3(p^s)\vee P^3(p^t)\ar[d]\ar[r]									&C\ar[d]^-{d}\\
\ast\ar[r]										&P^3(p^s)\times P^3(p^t)\ar@{=}[r]									&P^3(p^s)\times P^3(p^t)
}
\end{equation}
where $C_f$ is the homotopy cofibre of $f$ and $d$ is an induced map. As $C_{f}$ is $2$-connected, 
there is an isomorphism
\[
d^{\ast}:H^2(P^3(p^s)\times P^3(p^t);\mathbb{Z}/p^{u}\mathbb{Z})\to H^2(C;\mathbb{Z}/p^{u}\mathbb{Z}).
\]
Therefore $H^2(C;\mathbb{Z}/p^{u}\mathbb{Z})$ is 
generated by~$d^{\ast}(u_s\otimes 1)$ and $d^{\ast}(1\otimes u_t)$. 

The right column of~(\ref{CD_compare C_f and P^3xP^3}) induces the exact sequence
\begin{equation} 
\label{Cfexact1}
H^3(C;\mathbb{Z}/p^{u}\mathbb{Z})\rightarrow H^3(C_f;\mathbb{Z}/p^{u}\mathbb{Z})\overset{b}{\rightarrow}H^4(P^3(p^s)\times P^3(p^t);\mathbb{Z}/p^{u}\mathbb{Z})\overset{d^{\ast}}{\rightarrow}H^4(C;\mathbb{Z}/p^{u}\mathbb{Z}).
\end{equation} 
By Lemma~\ref{H4lemma}, 
$H^{4}(P^{3}(p^{s})\times P^{3}(p^{t});\mathbb{Z}/p^{u}\mathbb{Z})\cong\mathbb{Z}/p^{u}\mathbb{Z}$ 
is generated by the cup product~\mbox{$u_{s}\cup u_{t}$}. The naturality of the cup product implies 
that $d^{\ast}(u_{s}\cup u_{t})=d^{\ast}(u_{s})\cup d^{\ast}(u_{t})$. But by assumption, 
cup products in $\widetilde{H}^{\ast}(C;\mathbb{Z}/p^{u}\mathbb{Z})$ are zero. Therefore $d^{\ast}=0$ 
in~(\ref{Cfexact1}), implying that $b$ is onto. Hence the order of $H^{3}(C_{f};\mathbb{Z}/p^{u}\mathbb{Z})$ 
is at least $p^{u}$. 

On the other hand, the top row of~(\ref{CD_compare C_f and P^3xP^3}) induces the exact sequence
\begin{equation} 
\label{Cfexact2}
H^2(S^3;\mathbb{Z}/p^{u}\mathbb{Z})\rightarrow H^3(C_f;\mathbb{Z}/p^{u}\mathbb{Z})\overset{a}{\rightarrow}H^3(\Sigma P^2(p^s)\wedge P^2(p^t);\mathbb{Z}/p^{u}\mathbb{Z})\overset{f^*}{\rightarrow}H^3(S^3;\mathbb{Z}/p^{u}\mathbb{Z}).
\end{equation} 
Since $H^{2}(S^{3};\mathbb{Z}/p^{u}\mathbb{Z})=0$, the map $a$ is an injection, and by 
Lemma~\ref{Mooreskel}, 
\[H^{3}(\Sigma P^{2}(p^{s})\wedge P^{2}(p^{t});\mathbb{Z}/p^{u}\mathbb{Z})\cong\mathbb{Z}/p^{u}\mathbb{Z}.\] 
Hence the order of $H^{3}(C_{f};\mathbb{Z}/p^{u}\mathbb{Z})$ is at most $p^{u}$. 

Thus $H^{3}(C_{f};\mathbb{Z}/p^{u}\mathbb{Z})$ has order $p^{u}$. But this implies that $a$ 
is a monomorphism between finite groups of the same order and 
so must be an isomorphism. Therefore $f^{\ast}$ in~(\ref{Cfexact2}) is the zero map.  
\end{proof} 

A similar argument to Proposition~\ref{S3cuptriv}, but with variations, gives the following. 

\begin{prop} 
\label{Moorecuptriv}
Let $p$ be an odd prime and $r,s,t\geq 1$. Let  
$f:P^{3}(p^{r})\to\Sigma P^2(p^s)\wedge P^2(p^t)$  
be a map and let $C$ be the homotopy cofiber of the composite
\[
P^{3}(p^{r})\overset{f}{\longrightarrow}\Sigma P^2(p^s)\wedge P^2(p^t)
    \overset{[\imath_1,\imath_2]}{\longrightarrow}P^3(p^s)\vee P^3(p^t).
\]
Let $v=\min(r,s,t)$. Then the following are equivalent:
\begin{itemize} 
   \item[(a)] the map $f$ is null homotopic;
   \item[(b)] \(\namedright{H^{3}(\Sigma P^{2}(p^{s})\wedge P^{2}(p^{t});\mathbb{Z}/p^{v}\mathbb{Z})} 
                        {f^{\ast}}{H^{3}(P^{3}(p^{r});\mathbb{Z}/p^{v}\mathbb{Z})}\) 
                  is the zero map;
   \item[(c)] all cup products in $\widetilde{H}^*(C;\Z/p^{v}\Z)$ are zero.

\end{itemize} 
\end{prop}

\begin{proof} 
$(a)\Leftrightarrow(b)$: Let $u=\min(s,t)$ and consider the following string of isomorphisms 
\[\begin{split} 
      [P^{3}(p^{r}),\Sigma P^{2}(p^{s})\wedge P^{2}(p^{t})] 
           & \cong H_{3}(\Sigma P^{2}(p^{s})\wedge P^{2}(p^{t});\mathbb{Z}/p^{r}\mathbb{Z}) \\ 
       & \cong H^{3}(\Sigma P^{2}(p^{s})\wedge P^{2}(p^{t});\mathbb{Z}/p^{r}\mathbb{Z}) \\ 
       & \cong H^{3}(P^{4}(p^{u})\vee P^{5}(p^{u});\mathbb{Z}/p^{r}\mathbb{Z}) \\ 
       & \cong\left\{\begin{array}{ll} 
              \mathbb{Z}/p^{r}\mathbb{Z} & \mbox{if $r<u$} \\ 
              \mathbb{Z}/p^{u}\mathbb{Z} & \mbox{if $r\geq u$}\end{array}\right. \\ 
        & \cong\mathbb{Z}/p^{v}\mathbb{Z} \\ 
        & \cong H^{3}(P^{4}(p^{u})\vee P^{5}(p^{u});\mathbb{Z}/p^{v}\mathbb{Z}) \\ 
        & \cong H^{3}(\Sigma P^{2}(p^{r}\wedge P^{2}(p^{s});\mathbb{Z}/p^{v}\mathbb{Z}) 
  \end{split}\] 
The first isomorphism is due to the mod-$p^{r}$ Hurewicz isomorphism since 
$\Sigma P^{2}(p^{s})\wedge P^{2}(p^{t})$ is~$2$-connected. The second isomorphism 
holds by the Universal Coefficient Theorem and the third holds by Lemma~\ref{Mooreskel}. 
The fourth isomorphism is the calculation of degree~$3$ cohomology, the fifth holds 
since $v=\min(r,s,t)=\min(r,u)$, the sixth is calculation again, and the seventh holds 
by Lemma~\ref{Mooreskel}. The transition from the second to the seventh is induced 
by the map of coefficient rings induced by the epimorphism 
\(\namedright{\mathbb{Z}/p^{r}\mathbb{Z}}{}{\mathbb{Z}/p^{v}\mathbb{Z}}\). 
Thus, under these isomorphisms, a map 
\(f\colon\namedright{P^{3}(p^{r})}{}{\Sigma P^{2}(p^{s})\wedge P(p^{t})}\) 
is sent to the map it induces in mod-$p^{v}$ cohomology. Thus $f$ is null homotopic 
if and only if $f^{\ast}=0$ in mod-$p^{v}$ cohomology. 
\medskip 

$(b)\Leftrightarrow(c)$:
Consider the homotopy cofibration diagram  
\[
\xymatrix{
P^3(p^r)\ar@{=}[d]\ar[r]^-{f}						&\Sigma P^2(p^s)\wedge P^2(p^t)\ar[d]^-{[\imath_1,\imath_2]}\ar[r]	&C_f\ar[d]\\
P^3(p^r)\ar[d]\ar[r]^-{[\imath_1,\imath_2]\circ f}	&P^3(p^s)\vee P^3(p^t)\ar[d]\ar[r]									&C\ar[d]^-{d}\\
\ast\ar[r]										&P^3(p^s)\times P^3(p^t)\ar@{=}[r]									&P^3(p^s)\times P^3(p^t)
}
\] 
where $C_f$ is the homotopy cofibre of $f$ and $d$ is an induced map.
As $C_f$ is 2-connected,
\[
d^*:H^2(P^3(p^s)\times P^3(p^t);\Z/p^v\Z)\to H^2(C;\Z/p^v\Z)
\]
is an 
isomorphism. Therefore $H^2(C;\Z/p^v\Z)$ is generated by $d^*(u_s\otimes 1)$ and 
$d^*(1\otimes u_t)$. The diagram also induces a diagram of exact sequences 
\[
\xymatrix{
H^3(C_f;\Z/p^v\Z)\ar[r]^-{a}\ar[d]^-{b}	&H^3(\Sigma P^2(p^s)\wedge P^2(p^t);\Z/p^v\Z)\ar[r]^-{f^*}\ar[d]^-{c}	&H^3(P^3(p^r);\Z/p^v\Z)\\
H^4(P^3(p^s)\times P^3(p^t);\Z/p^v\Z)\ar@{=}[r]\ar[d]^-{d^{\ast}} &H^4(P^3(p^s)\times P^3(p^t);\Z/p^v\Z)\ar[d]	&\\
H^4(C;\Z/p^v\Z)\ar[r]	&H^4(P^3(p^s)\vee P^3(p^s);\Z/p^v\Z)=0	&
}
\] 
where $a$, $b$ and $c$ are names for the maps induced in cohomology. 
Observe that, in the middle column, $s,t\geq v$ so
\[
H^3(\Sigma P^2(p^s)\wedge P^2(p^t);\Z/p^v\Z)\cong H^4(P^3(p^s)\times P^3(p^t);\Z/p^v\Z)\cong\Z/p^v\Z,
\]
implying that $c$ is an isomorphism. Therefore, the commutativity of the top square implies that $a$ is surjective if and only if $b$ is. On the other hand, the top row implies that~$a$ is surjective if and only if $f^*$ is the zero map, while the left column implies that~$b$ is surjective if and only if $d^{\ast}$ is the zero map. Thus $f^{\ast}=0$ if and only in $d^{\ast}=0$. 
Since~\mbox{$H^4(P^3(p^s)\times P^3(p^t);\Z/p^v\Z)$} is generated by $u_s\cup u_t$, $d^{\ast}=0$ if and only if $H^4(C;\Z/p^v\Z)$ has no cup products. Hence $f^{\ast}=0$ if and only if $H^4(C;\Z/p^v\Z)$ has no cup products.
\end{proof}

\section{A homological test for a null homotopy II} 
\label{sec:test2} 

In this section we aim towards Proposition~\ref{trivtestprop}, which gives homological 
and cohomological criteria for when certain maps are null homotopic, and which is 
applicable much more widely than Propositions~\ref{S3cuptriv} and~\ref{Moorecuptriv}. 
It also generalizes a result of Neisendorfer~\cite[Corollary 11.12]{neisendorfer10} 
on the mod-$p^{r}$ Hopf invariant. We rephrase that result in weaker form for a better 
comparison to Proposition~\ref{trivtestprop}. 

\begin{lemma} 
   \label{modphopfinv} 
   Let $p$ be an odd prime and $r,s\geq 1$. Let 
   \(f\colon\namedright{P^{3}(p^{r})}{}{P^{3}(p^{s})}\) 
   be a map and let~$C_{f}$ be its cofibre. If 
   \begin{itemize} 
      \item \(f_{\ast}\colon\namedright{\tilde{H}_{\ast}(P^{3}(p^{r});\mathbb{Z})}{} 
                     {\tilde{H}_{\ast}(P^{3}(p^{s});\mathbb{Z})}\) 
               is the zero map, and 
      \item all cup products in $\widetilde{H}^{\ast}(C_{f};\mathbb{Z}/p^{min(r,s)}\mathbb{Z})$ are zero,
   \end{itemize} 
   then $f$ is null homotopic.~$\qqed$  
\end{lemma} 

Lemma~\ref{modphopfinv} will be generalized to maps 
\(f\colon\namedright{X}{}{\bigvee_{i=1}^{m} P^{3}(p^{r_i})}\) 
for $X=S^{3}$ or $X=P^{3}(p^{r})$. This requires some initial work, the 
first aspect of which is a general lemma concerning trivial cup products 
related to maps of wedges. 

\begin{lemma} 
\label{pinchcupabstract} 
   Let 
   \(f\colon\namedright{\bigvee_{i=1}^{m} A_{i}}{}{\bigvee_{j=1}^{n} B_{j}}\) 
   be a map with homotopy cofibre $C_{f}$ and suppose that $f^{\ast}=0$ 
   for cohomology with coefficient group $G$ and all cup products in $\widetilde{H}^{\ast}(C_{f};G)$ 
   are zero. For $1\leq\imath\leq m$ and $1\leq\jmath\leq n$, let $f_{\imath,\jmath}$ be 
   the composite 
   \[f_{\imath,\jmath}\colon A_{\imath}\hookrightarrow 
           \nameddright{\bigvee_{i=1}^{m} A_{i}}{f}{\bigvee_{j=1}^{n} B_{j}}{}{B_{\jmath}}\] 
   where the left map is the inclusion of the $\imath^{th}$ wedge summand and 
   the right map is the pinch onto the $\jmath^{th}$ wedge summand. If $C_{f_{\imath,\jmath}}$ 
   is the homotopy cofibre of $f_{\imath,\jmath}$ then all cup products in~$\widetilde{H}^{\ast}(C_{f_{\imath,\jmath}};G)$ are zero. 
\end{lemma}  

\begin{proof} 
We use an intermediate map. Let $f_{\jmath}$ be the composite 
\begin{equation} 
   \label{fjdef} 
   f_{\jmath}\colon\nameddright{\bigvee_{i=1}^{m} A_{i}}{f}{\bigvee_{j=1}^{n} B_{j}}{}{B_{\jmath}} 
\end{equation}  
and let $C_{f_{\jmath}}$ be the homotopy cofibre of $f_{\jmath}$. 
Consider the homotopy cofibration diagram
\[\diagram 
     \bigvee_{i=1}^{m} A_{i}\rto^-{f}\ddouble & \bigvee_{j=1}^{n} B_{j}\rto\dto & C_{f}\dto^{d} \\ 
     \bigvee_{i=1}^{m} A_{i}\rto^-{f_{\jmath}} & B_{\jmath}\rto & C_{f_{\jmath}}  
  \enddiagram\] 
where $d$ is an induced map of cofibres. Take cohomology with coefficient group $G$.  
The homotopy cofibration diagram induces a map between long exact sequences in 
cohomology. By hypothesis, $f^{\ast}=0$ so the definition of $f_{\jmath}$ implies 
that $f_{\jmath}^{\ast}=0$ as well. Therefore, for every~$k\geq 1$, there is a commutative 
diagram of exact sequences 
\[\diagram 
     0\rto\ddouble & H^{k}(\bigvee^m_{i=1}\Sigma A;G)\rto\ddouble & H^{k}(C_{f_{\jmath}}:G)\rto\dto^{d^{\ast}} 
          & H^{k}(B_{\jmath};G)\rto\dto & 0\ddouble \\ 
     0\rto & H^{k}(\bigvee^m_{i=1}\Sigma A;G)\rto & H^{k}(C_{f};G)\rto & H^{k}(\bigvee_{i=1}^{n} B_{i};G)\rto & 0. 
  \enddiagram\]
A diagram chase shows that $d^{\ast}$ is injective, and this is true for all $k\geq 1$. Thus, 
by the naturality of the cup product,  the vanishing of cup products in $\widetilde{H}^{\ast}(C_{f};G)$ 
implies their vanishing in $\widetilde{H}^{\ast}(C_{f_{\jmath}};G)$.

Next, notice that the definition of $f_{\imath,\jmath}$ in the statement of the lemma 
and $f_{\jmath}$ in~(\ref{fjdef}) imply that $f_{\imath,\jmath}$ is the composite 
\(A_{\imath}\hookrightarrow\namedright{\bigvee_{i=1}^{m} A_{i}}{f_{\jmath}}{B_{\jmath}}\). 
This factorization induces a homotopy cofibration diagram  
\[\diagram 
      A_{\imath}\rto\ddouble & \bigvee_{i=1}^{m} A_{i}\rto^{h}\dto^{f_{\jmath}} 
              & \bigvee_{\substack{i=1\\ i\neq\imath}}^{m} A_{i}\dto^{g} \\ 
      A_{\imath}\rto^-{f_{\imath,\jmath}}\dto & B_{\jmath}\rto\dto & C_{f_{\imath,\jmath}}\dto^{d'} \\ 
\ast\rto      & C_{f_{\jmath}}\rdouble & C_{f_{\jmath}}  
  \enddiagram\]
where $h$ is the pinch map, and $g$ and $d'$ are induced maps. Since $f_{\jmath}^{\ast}=0$ and
\[
h^*\colon H^*(\bigvee_{\substack{i=1\\ i\neq\imath}}^{m}A_i;G)\to H^*(\bigvee^m_{i=1}A_i;G)
\]
is an injection, the top right square implies that $g^{\ast}=0$. Therefore, from the right vertical cofibration in the 
preceding diagram we obtain a surjection 
\[
(d')^*\colon H^{\ast}(C_{f_{\jmath}};G)\to H^{\ast}(C_{f_{\imath,\jmath}};G).
\] 
As cup products in $\widetilde{H}^{\ast}(C_{f_{\jmath}};G)$ are zero and $(d')^{\ast}$ is a surjection, 
cup products in~$\widetilde{H}^{\ast}(C_{f_{\imath,\jmath}};G)$ are also zero. 
\end{proof}

Next, we make a transition from a hypothesis that a map is zero in cohomology as in 
Lemma~\ref{pinchcupabstract} to a map being zero in homology. In general, 
if the coefficient group $G$ in Lemma~\ref{pinchcupabstract} is a field then 
the Universal Coefficient Theorem immediately implies that if~\mbox{$f_{\ast}=0$} then $f^{\ast}=0$. 
The coefficient ring we care about is $\mathbb{Z}/p^{r}\mathbb{Z}$, so we need 
to be more cautious. Perhaps overdoing it, we focus on the $3$-dimensional Moore 
space case again. 

\begin{lemma}\label{lemma_f_* trivial implies all f^* trivial}
Let $p$ be an odd prime and let $r\geq 1$. Let $X=P^3(p^r)$ or $S^3$ and let 
$f:X\to\bigvee^m_{i=1}P^3(p^{r_i})$ 
be a map. If $f_*:\tilde{H}_*(X;\Z)\to\tilde{H}_*(\bigvee_{i=1}^mP^3(p^{r_i});\Z)$ is trivial then 
for any abelian group $G$ the map $f^*:\tilde{H}^*(\bigvee_{i=1}^mP^3(p^{r_i});G)\to\tilde{H}^*(X;G)$ 
is trivial.
\end{lemma} 

\begin{proof}
It suffices to prove the lemma in the $m=1$ case. 
For $X=P^3(p^r)$ it is obvious that $f^*:\tilde{H}^j(P^3(p^{r_1});G)\to\tilde{H}^j(P^3(p^r);G)$ 
is trivial except possibly for $j\in\{2,3\}$. By the Universal Coefficient Theorem,  
there are natural isomorphisms 
\[
H^2(P^3(p^r);G)\cong Hom(H_2(P^3(p^r);\mathbb{Z}),G))
\]
and
\[
H^3(P^3(p^r);G)\cong Ext(H_2(P^3(p^r);\mathbb{Z}),G)).
\] 
By hypothesis, 
$f_{\ast}\colon H_{2}(P^{3}(p^{r});\mathbb{Z})\to H_{2}(P^{3}(p^{r_1});\mathbb{Z})$  
is the zero map, so the naturality of the Universal Coefficient Theorem implies that  
$f^{\ast}\colon H^{j}(P^{3}(p^{r_1});G)\to H^{j}(P^{3}(p^{r});G)$ is the zero map 
for \mbox{$j\in\{2,3\}$}. 

For $X=S^3$, it suffices to show that $f^*:H^3(P^3(p^{r_1});G)\to H^3(S^3;G)$ is trivial. 
Let $\rho:P^3(p^{r_1})\to S^3$ be the pinch map to the top cell and consider the composite 
\begin{equation} 
  \label{rhof}  
  H^{3}(S^{3};G)\overset{\rho^{\ast}}{\longrightarrow} H^{3}(P^{3}(p^{r_1});G)\overset{f^{\ast}}{\longrightarrow} 
      H^{3}(S^{3};G). 
\end{equation}  
Observe that the long exact sequence in cohomology determined by the homotopy cofibration~
$S^{2}\to P^{3}(p^{r_1})\overset{\rho}{\longrightarrow} S^{3}$ 
implies that $\rho^{\ast}$ in~(\ref{rhof}) is an epimorphism. Therefore, in~(\ref{rhof}), $f^{\ast}=0$ 
if and only if $f^{\ast}\circ\rho^{\ast}=0$. But $\rho\circ f$ is a self-map of $S^{3}$ which factors 
through a rationally contractible space, implying that it is null homotopic. Hence $f^{\ast}\circ\rho^{\ast}=0$, 
and so $f^{\ast}=0$. 
\end{proof} 

In general, the Hilton-Milnor Theorem states that there is a homotopy equivalence 
\begin{equation} 
  \label{HM} 
  \Omega(\bigvee_{i=1}^{m}\Sigma Y_{i})\simeq\prod_{\alpha\in\mathcal{I}} 
       \Omega\Sigma(Y_{1}^{\wedge\alpha_{1}}\wedge\cdots\wedge Y_{m}^{\wedge\alpha_{m}}) 
\end{equation} 
where $\mathcal{I}$ runs over a module basis for the free Lie algebra $L\langle v_{1},\ldots,v_{m}\rangle$, 
and if $\alpha\in L\langle v_{1},\ldots,v_{m}\rangle$ is a module basis element then for $1\leq i\leq m$ 
the integer $\alpha_{i}$ records the number of instances of $v_{i}$ in $\alpha$. Here, if $\alpha_{i}=0$ 
for some $i$ then the smash product 
$Y_{1}^{\wedge\alpha_{1}}\wedge\cdots\wedge Y_{m}^{\wedge\alpha_{m}}$ 
is regarded as omitting $Y_{i}$ rather than being a point; for example, 
$Y_{1}^{\wedge 2}\wedge Y_{2}^{\wedge 0}\wedge Y_{3}^{\wedge 3}$ is regarded 
as $Y_{1}^{\wedge 2}\wedge Y_{3}^{\wedge 3}$. Moreover, for $1\leq k\leq m$ let 
\[\iota_{k}\colon\namedright{\Sigma Y_{k}}{}{\bigvee_{i=1}^{m}\Sigma Y_{i}}\] 
be the inclusion of the $k^{th}$ wedge summand. For $\alpha\in\mathcal{I}$, let  
\[w_{\alpha}\colon\namedright{\Sigma(Y_{1}^{\wedge\alpha_{1}}\wedge\cdots\wedge Y_{m}^{\wedge\alpha_{m}})} 
      {}{\bigvee_{i=1}^{m}\Sigma Y_{i}}\] 
be the iterated Whitehead product formed from the maps $\iota_{k}$ where each 
instance of $v_{k}$ in $\alpha$ is represented by the map $\iota_{k}$. Then the homotopy 
equivalence~(\ref{HM}) is realized by multiplying together the maps $\Omega w_{\alpha}$ 
using the loop structure on $\Omega(\bigvee_{i=1}^{m}\Sigma Y_{i})$. 

In our case, we have 
\[\Omega(\bigvee_{i=1}^{m} P^{3}(p^{r_{i}}))\simeq\prod_{\alpha\in\mathcal{I}} 
       \Omega\Sigma P^{2}(p^{r_1})^{\wedge\alpha_{1}}\wedge\cdots\wedge P^{2}(p^{r_m})^{\wedge\alpha_{m}}.\] 
Observe that $P^{2}(p^{r_1})^{\wedge\alpha_{1}}\wedge\cdots\wedge P^{2}(p^{r_m})^{\wedge\alpha_{m}}$ 
is $((\alpha_{1}+\cdots+\alpha_{m})-1)$-connected. Suppose that~$X'$ is $2$-dimensional. Then 
$[X',\Omega\Sigma P^{2}(p^{r_1})^{\wedge\alpha_{1}}\wedge\cdots\wedge P^{2}(p^{r_m})^{\wedge\alpha_{m}}]\cong 0$ 
if $(\alpha_{1}+\cdots+\alpha_{m})\geq 3$. Observe also that there are $m$ cases 
for which $(\alpha_{1}+\cdots +\alpha_{m})=1$ and $\binom{m}{2}$ cases for which~$(\alpha_{1}+\cdots +\alpha_{m})=2$. So if $X=\Sigma X'$ then 
\begin{eqnarray*}
[X, \bigvee_{i=1}^{m} P^3(p^{r_i})]
  &\cong 
    &[X', \Omega(\bigvee_{i=1}^{m} P^3(p^{r_i}))]\\
  &\cong 
    &[X', \prod_{j=1}^{m}\Omega P^3(p^{r_j})\times 
          \prod_{k\neq l}\Omega\Sigma P^2(p^{r_k})\wedge P^2(p^{r_l})]\\
  &\cong &\prod_{j=1}^{m} [X, P^3(p^{r_j})]\times\prod_{k\neq l}[X,\Sigma P^2(p^{r_k})\wedge P^2(p^{r_l})].  
\end{eqnarray*}
Further, the $j^{th}$ factor $[X,P^3(p^{r_j})]$ is mapped to $[X, \bigvee_{i=1}^{m} P^3(p^{r_i})]$ by 
the inclusion $\iota_{j}$ and the~$\binom{m}{2}$ factors $[X,\Sigma P^2(p^{r_k})\wedge P^2(p^{r_l})]$
may be arranged so that they map to $[X, \bigvee_{i=1}^{m} P^3(p^{r_i})]$ by the Whitehead products
\[
\Sigma P^2(p^{r_k})\wedge P^2(p^{r_l})\overset{[\iota_{k},\iota_{l}]}{\longrightarrow}P^3(p^{r_k})\vee P^3(p^{r_l})\hookrightarrow\bigvee_{i=1}^{m} P^3(p^{r_i}) 
\]
where $1\leq k<l\leq m$. Hence if 
\(f\colon\namedright{X}{}{\bigvee_{i=1}^{m} P^{3}(p^{r_i})}\) 
then we may write 
\begin{equation} 
  \label{fdecomp}
  f\simeq\sum_{j=1}^{m}\iota_{j}\circ g_j+\sum_{1\leq k< l\leq m}[\iota_{k},\iota_{l}]\circ h_{k,l}
\end{equation}   
for maps
\(\namedright{X}{g_{j}}{}{P^{3}(p^{r_j})}\)
and 
\(\namedright{X}{h_{k, l}}{\Sigma P^{2}(p^{r_k})\wedge P^{2}(p^{r_l})}\). 

\begin{prop} 
   \label{trivtestprop}
   Let $X=P^3(p^r)$ where $p$ is an odd prime and $r\geq 1$ or let $X=S^3$ and set $r=\infty$. Let
   $f\colon X\to\bigvee_{i=1}^{m} P^3(p^{r_i})$ be a map and let $C_f$ be its cofiber. If 
   \begin{itemize} 
      \item $f_*:\tilde{H}_*(X;\Z)\to\tilde{H}_*(\bigvee_{i=1}^{m} P^3(p^{r_i});\Z)$ is 
               the zero map and 
      \item all cup products in $\widetilde{H}^*(C_f;\Z/p^{min(r,r_i)}\Z)$ are zero for all $1\leq i\leq m$,
      \end{itemize} 
   then $f$ is null homotopic. 
\end{prop}

\begin{proof} 
Since $X$ is $S^{3}$ or $P^{3}(p^{r})$ we have $X\simeq\Sigma X'$ where 
$X'$ is $2$-dimensional. Therefore, by~(\ref{fdecomp}), we have  
$f\simeq\sum_{j=1}^{m}\iota_{j}\circ g_j+\sum_{1\leq k< l\leq m}[\iota_{k},\iota_{l}]\circ h_{k,l}$  
for maps
\(\namedright{X}{g_{j}}{}{P^{3}(p^{r_j})}\) 
and 
\(\namedright{X}{h_{k, l}}{\Sigma P^{2}(p^{r_k})\wedge P^{2}(p^{r_l})}\). 
To show that $f$ is null homotopic it suffices to show that each $g_j$ and $h_{k, l}$ 
is null homotopic. 

First consider the map $g_{j}$ when $X=P^{3}(p^{r})$. Notice that $g_j$ is the composite
\[
g_j:P^{3}(p^{r})\overset{f}{\longrightarrow}\bigvee_{i=1}^{m} P^3(p^{r_i})\overset{q}{\longrightarrow}P^3(p^{r_j})
\]
where $q$ is the pinch map onto the $j^{th}$ wedge summand. 
Since $f$ induces the zero map in integral homology, so does $g_{j}$. The spaces 
involved let us apply Lemma~\ref{lemma_f_* trivial implies all f^* trivial}, showing that ~$g_{j}$ induces the zero map in mod-$p^{min(r,r_j)}$ cohomology. By hypothesis, 
all cup products in $\widetilde{H}^{\ast}(C_{f};\mathbb{Z}/p^{min(r,r_j)}\mathbb{Z})$ are zero, so by 
Lemma~\ref{pinchcupabstract}, all cup products in $\widetilde{H}^{\ast}(C_{g_{j}};\mathbb{Z}/p^{min(r,r_j)}\mathbb{Z})$ 
are also zero. Thus, by Lemma~\ref{modphopfinv}, $g_{j}$ is null homotopic. 

Next, consider the map $g_{j}$ when $X=S^{3}$. Now $g_{j}$ is the composite
\(\nameddright{S^{3}}{f}{\bigvee_{i=1}^{m} P^{3}(p^{r_i})}{q}{P^{3}(p^{r_j})}\). 
Consider the composite 
\[\overline{g}_{j}\colon\nameddright{P^{3}(p^{r_j})}{\pi}{S^{3}}{g_{j}}{P^{3}(p^{r_j})}\]
where $\pi$ is the pinch map to the top cell. The argument in the previous paragraph 
implies that~$\overline{g}_{j}$ is null homotopic. Therefore $g_{k}$ extends across 
the cofibre of $\pi$, implying that $g_{k}$ factors as a composite 
\(\nameddright{S^{3}}{p^{r_j}}{S^{3}}{\gamma_{j}}{P^{3}(p^{r_j})}\) 
for some map $\gamma_{j}$. By Lemma~\ref{pi3Moore}, 
$\pi_{3}(P^{3}(p^{r_j}))\cong\mathbb{Z}/p^{r_j}\mathbb{Z}$, 
so $g_{j}\simeq p^{r_j}\cdot\gamma_{j}$ is null homotopic. 

At this point, we have shown that for either $X=S^{3}$ or $P^{3}(p^{r})$ we have 
$g_{j}$ null homotopic for $1\leq j\leq m$. Thus~(\ref{fdecomp}) implies that 
$f\simeq\sum_{1\leq k< l\leq m}[\iota_{k},\iota_{l}]\circ h_{k,l}$. Let 
\[q_{k,l}\colon\namedright{\bigvee_{i=1}^{m} P^{3}(p^{r_i})}{}{P^{3}(p^{r_k})\vee P^{3}(p^{r_l})}\] 
be the pinch map onto the $k^{th}$ and $l^{th}$ wedge summands. Observe that 
every Whitehead product $[\iota_{s},\iota_{t}]$ for $1\leq s<t\leq m$ composes trivially 
with $q_{k, l}$ except $[\iota_{k},\iota_{l}]$. Therefore~\mbox{$q_{k, l}\circ f\simeq q_{k, l}\circ(\sum_{1\leq s< t\leq m}[\iota_{s},\iota_{t}]\circ h_{s,t})\simeq 
       [\iota_{k},\iota_{l}]\circ h_{k, l}$}. 
That is, $q_{k,l}\circ f$ is homotopic to the composite 
\[\overline{h}_{k,l}\colon\nameddright{X}{h_{k,l}}{\Sigma P^{2}(p^{r_k})\wedge P^{2}(p^{r_l})} 
        {[\iota_{k},\iota_{l}]}{P^{3}(p^{r_k})\vee P^{3}(p^{r_l})}.\]
Since $f$ induces the zero map in integral homology, so does $\overline{h}_{k,l}$. 
Let $C_{\overline{h}_{k,l}}$ be the homotopy cofibre of $\overline{h}_{k,l}$. 
By hypothesis, cup products in $\widetilde{H}^{\ast}(C_{f};\mathbb{Z}/p^{min(r,r_i)}\mathbb{Z})$ are zero 
 for~$1\leq i\leq m$ so cup products in $\widetilde{H}^{\ast}(C_{f};\mathbb{Z}/p^{min(r,r_k,r_l)}\mathbb{Z})$ are zero. By 
Lemma~\ref{pinchcupabstract} (with~\mbox{$B_{\jmath}=P^{3}(p^{r_k})\vee P^{3}(p^{r_l})$}), cup products 
in $\widetilde{H}^{\ast}(C_{\overline{h}_{k,l}};\mathbb{Z}/p^{min(r,r_k,r_l)}\mathbb{Z})$
are also zero. Therefore, by Proposition~\ref{S3cuptriv} in the case $X=S^{3}$ and 
Proposition~\ref{Moorecuptriv} in the case $X=P^{3}(p^{r})$, the map~$h_{k,l}$ is 
null homotopic. As this is true for all $1\leq k<l\leq m$ we obtain $f\simeq\ast$. 
\end{proof}

\section{The homotopy type of the suspension of certain $CW$-complexes} 
\label{sec:suspM} 

In this section we assume $M$ to be a 4-dimensional finite $CW$-complex that has one 4-cell and homology as follows: 
\begin{equation} 
\label{Mhlgy} 
\begin{array}{c|c}
i	&H_i(M;\mathbb{Z})\\
\hline
0	&\Z\\[1mm]
1	&\Z^\ell\oplus\bigoplus^n_{j=1}\Z/b_j\Z\\[2mm]
2	&\Z^d\oplus\bigoplus^{\bar{n}}_{\bar{j}=1}\Z/\bar{b}_{\bar{j}}\Z\\[2mm]
3	&\hspace{3mm}\Z^m\\[1mm]
4	&\Z\\[1mm]
\geq5	&0
\end{array}
\end{equation} 
Here each $b_{j}$ and $\bar{b}_{\bar{j}}$ is a power of an odd prime.

First consider the integer summands of $H_{1}(M;\mathbb{Z})$. 
Since the Hurewicz homomorphism $\pi_1(M)\to H_1(M;\mathbb{Z})$ is an epimorphism, 
each direct summand $\Z$ of $H_1(M;\mathbb{Z})$ is generated by the Hurewicz image 
of some map 
\(\alpha_{i}\colon\namedright{S^{1}}{}{M}\).
Let 
\[a\colon\namedright{\bigvee_{i=1}^{\ell} S^{1}}{}{M}\] 
be the wedge sum of the maps $\alpha_{i}$ and let $W$ be the homotopy cofibre of $a$. 
  
\begin{lemma} 
   \label{S1M} 
   The map $\Sigma a$ has a left homotopy inverse and there is a homotopy equivalence 
   \[\Sigma M\simeq(\bigvee_{i=1}^{\ell} S^{2})\vee\Sigma W.\] 
\end{lemma} 

\begin{proof} 
The Hurewicz Theorem implies that the image of $a_{\ast}$ 
is $H_{1}(M;\mathbb{Z})_{free}\cong\mathbb{Z}^{\ell}$. The Universal Coefficient 
Theorem implies that $H^{1}(M;\mathbb{Z})_{free}\cong H_{1}(M;\mathbb{Z})_{free}$. 
Let $a_{i}\in H_{1}(M;\mathbb{Z})$ be the image of $(\alpha_{i})_{\ast}$ 
and $\bar{a}_{i}\in H^{1}(M;\mathbb{Z})$ be the dual of $a_{i}$. Then $\bar{a}_{i}$ is 
represented by a map 
\(\epsilon_{i}\colon\namedright{M}{}{K(\mathbb{Z},1)\simeq S^{1}}\) 
and the composite 
\(\nameddright{S^{1}}{\alpha_{i}}{M}{\epsilon_{i}}{S^{1}}\) 
is the identity map. After suspending one may use the co-H structure to give a map
\(\epsilon\colon\namedright{\Sigma M}{}{\bigvee_{i=1}^{\ell} S^{2}}\) 
which is a left homotopy inverse for $\Sigma a$. Therefore, with respect to the 
homotopy cofibration,~\(\nameddright{\bigvee_{i=1}^{\ell} S^{2}}{\Sigma a}{\Sigma M}{\Sigma w}{\Sigma W}\) where $w:M\to W$ is the quotient map, 
if $\sigma$ is the comultiplication on $\Sigma M$, the composite 
\[e\colon\lnameddright{\Sigma M}{\sigma}{\Sigma M\vee\Sigma M}{\epsilon\vee\Sigma w} 
      {(\bigvee_{i=1}^{\ell} S^{2})\vee\Sigma W}\] 
induces an isomorphism in homology. As the domain and range of $e$ are simply-connected, 
Whitehead's Theorem implies that $e$ is a homotopy equivalence.  
\end{proof} 

The description of $H_{\ast}(M;\mathbb{Z})$ in~(\ref{Mhlgy}) implies that the homology of $W$ is as follows: 
\[
\begin{array}{c|c}
i	&H_i(W;\mathbb{Z})\\
\hline
0	&\Z\\[1mm]
1	&\bigoplus^n_{j=1}\Z/b_j\Z\\[2mm]
2	&\Z^d\oplus\bigoplus^{\bar{n}}_{\bar{j}=1}\Z/\bar{b}_{\bar{j}}\Z\\[2mm]
3	&\hspace{3mm}\Z^m\\[1mm]
4	&\Z\\[1mm]
\geq5	&0
\end{array}
\]
We wish to give a homotopy decomposition of $\Sigma W$ as a wedge 
of spheres and Moore spaces. To do so we analyze the homology decomposition 
of $\Sigma W$.  

Define $M(\Z/k\Z, n)=P^{n+1}(k)$ and $M(\Z, n)=S^n$, and for any finitely generated 
abelian groups $A$ and $B$ define $M(A\oplus B, n)=M(A, n)\vee M(B, n)$. Then 
$\tilde{H}_i(M(A, n);\mathbb{Z})$ is $A$ for~$i=n$ and zero otherwise. The following 
lemma describes the homology decomposition of a simply-connected CW-complex. 

\begin{lemma}[Theorem 4H.3, \cite{hatcher}]\label{lemma_minimal cell structure}
Let $X$ be an $n$-dimensional simply-connected CW-complex and let $H_{i}=H_{i}(X;\mathbb{Z})$. 
Then there is a sequence of subcomplexes $\{X_i\}^n_{i=1}$ such that
\begin{enumerate}
\item	$H_i(X_m;\mathbb{Z})\cong H_i(X;\mathbb{Z})$ for $i\leq m$ and 
         $H_i(X_m;\mathbb{Z})=0$ for $i>m$;
\item	$X_2=M(H_2,2)$ and $X\simeq X_n$;
\item	$X_{m+1}$ is the mapping cone of a map $f_m\colon M(H_{m+1},m)\to X_m$  
         that induces a trivial homomorphism 
         $(f_m)_*\colon H_m(M(H_{m+1},m);\mathbb{Z})\to H_m(X_m;\mathbb{Z})$.
\end{enumerate} 
\end{lemma}

In our case, to describe the homology decomposition of $\Sigma W$ we need 
some notation. Let 
\[P=\bigvee_{j=1}^{n} P^{3}(b_{j})\qquad
    \overline{P}=\bigvee_{\bar{j}=1}^{\bar{n}} P^{3}(\bar{b}_{\bar{j}})\qquad 
    \mbox{and}\qquad S=\bigvee_{k=1}^{d} S^{2}.\] 
Starting with $W_{2}=P$, Lemma~\ref{lemma_minimal cell structure} implies that there 
are homotopy cofibrations\smallskip  
\begin{equation} 
  \label{hlgydecomp} 
  \begin{split} 
       \nameddright{S\vee\overline{P}}{f_{2}}{W_{2}}{}{W_{3}} \\ 
       \nameddright{\bigvee_{i=1}^{m} S^{3}}{f_{3}}{W_{3}}{}{W_{4}} \\ 
       \nameddright{S^{4}}{f_{4}}{W_{4}}{}{\Sigma W} 
  \end{split} 
\end{equation} 
where $f_{2}$, $f_{3}$ and $f_{4}$ induce the zero map in integral homology. 
In Lemmas~\ref{W3decomp} and \ref{W4decomp} we will show that the maps $f_{2}$ 
and $f_{3}$ are null homotopic, and in Lemma~\ref{suspWdecomp} we will show that the map~$f_{4}$ 
is either null homotopic or factors is an entirely controllable way. As this 
will involve analyzing maps between Moore spaces of different torsion orders, 
a preliminary lemma is required. 

\begin{lemma}\label{lemma_[Pa,Pb]=0}
Let $X$ be a finite CW-complex. If $p$ and $q$ are distinct primes and $m,n\geq 3$, then any map $f\colon P^m(p^r)\to\Sigma X\vee P^n(q^t)$ is homotopic to the composite
\[
P^m(p^r)\xrightarrow{f'}\Sigma X\hookrightarrow\Sigma X\vee P^n(q^t)
\]
where $f'$ is the composite 
\(\lnameddright{P^{m}(p^{r})}{f}{\Sigma C\vee P^{n}(q^{t})}{\mbox{\tiny\rm pinch}}{\Sigma X}\). 
\end{lemma}

\begin{proof}
First we show that $[P^m(p^r),Z\wedge P^m(q^t)]$ is trivial for any finite path-connected 
CW-complex $Z$. By the K\"{u}nneth Theorem there is an exact sequence
\[
0\to\bigoplus^n_{i=1}\tilde{H}_i(Z)\otimes\tilde{H}_{n-i}(P^n(q^t))\to\tilde{H}_n(Z\wedge P^n(q^t))\to\bigoplus^n_{i=1}\text{Tor}(\tilde{H}_i(Z),\tilde{H}_{n-i-1}(P^n(q^t)))\to0.
\]
This implies that the groups $\tilde{H}_*(Z\wedge P^n(q^t))$ are finite abelian and consist only of $q$-torsion. 
Therefore, by Serre's Theorem, the homotopy groups $\pi_i(Z\wedge P^n(q^t))$ are also finite abelian 
and consist only of $q$-torsion. The homotopy cofibration
\[
S^{m-1}\overset{p^r}{\longrightarrow}S^{m-1}\longrightarrow P^m(p^r)
\]
induces an exact sequence
\[
\pi_m(Z\wedge P^n(q^{t}))\overset{p^r}{\longrightarrow}\pi_m(Z\wedge P^n(q^{t}))\longrightarrow[P^m(p^{r}),Z\wedge P^n(q^{t})]  
       \longrightarrow
\] 
\[
\hspace{6cm} \pi_{m-1}(Z\wedge P^n(q^{t}))\overset{p^r}{\longrightarrow}\pi_{m-1}(Z\wedge P^n(q^{t})). 
\]
Since multiplying $\pi_i(Z\wedge P^n(q^t))$ by $p^r$ is an isomorphism for $i\geq 1$, 
by exactness we obtain $[P^m(p^r), Z\wedge P^n(q^t)]\cong0$.

Next, the homotopy class of $f$ is in $[P^m(p^r),\Sigma X\vee P^n(q^t)]$. Noting 
that both $P^{m}(p^{r})$ and $P^{n}(q^{t})$ are suspensions since $m,n\geq 3$, the 
Hilton-Milnor Theorem implies that 
\[
[P^m(p^r),\Sigma X\vee P^n(q^t)]\cong 
  \prod_{\alpha\in\mathcal{I}} [P^m(p^r),\Sigma X^{\wedge\alpha_{1}}\wedge (P^{n-1}(q^t))^{\wedge\alpha_{2}}]
\] 
where $\mathcal{I}$ runs over a Hall basis for the free Lie algebra $L\langle u,v\rangle$ and 
$\alpha_{1},\alpha_{2}$ count the number of instances of $u,v$ respectively in the bracket 
corresponding to $\alpha$. The argument in the first paragraph implies that if $\alpha_{2}\geq 1$ then 
each factor $[P^m(p^r),\Sigma X^{\wedge\alpha_{1}}\wedge (P^{n-1}(q^t))^{\wedge\alpha_{2}}]$, 
which is isomorphic to $[P^m(p^r), Z\wedge P^{n}(q^t)]$  
for $Z=X^{\wedge\alpha_{1}}\wedge (P^{n-1}(q^t))^{\wedge\alpha_{2}-1}$, 
equals zero.  The Hall basis for $L\langle u,v\rangle$ only has one term with $\alpha_{2}=0$, 
and that is $u$ (when $\alpha_{1}=1$). Thus 
\[
[P^m(p^r),\Sigma X\vee P^n(q^t)]\cong [P^m(p^r),\Sigma X].  
\] 
Hence $f$ factors through $f'$ up to homotopy.
\end{proof}

We also need a lemma concerning cup products in $W_{3}$. 

%a more general lemma concerning cup products associated to the homology decomposition 
%of a finite simply-connected $CW$-complex. 

%\begin{lemma} 
%   \label{hlgydecompcup} 
%   Let $Y$ be an $n$-dimensional $CW$-complex and let $\{Y_{i}\}_{i=2}^{n+1}$ 
%   be the homology decomposition of $\Sigma Y$. Then each of the maps 
%   \(\namedright{Y_{i}}{}{Y_{n+1}\simeq\Sigma Y}\) 
%   induces an epimorphism in cohomology. Consequently, cup products vanish in 
%   $H^{\ast}(Y_{i};\mathbb{Z})$ for all $2\leq i\leq n+1$. 
%\end{lemma} 

%\begin{proof} 
%The Universal Coefficient Theorem implies that, for all $m\geq 1$, there is an isomorphism 
%\[H^{m}(Y_{i};\mathbb{Z})\cong H_{m}(Y_{i};\mathbb{Z})_{free}\oplus H_{m-1}(Y_{i};\mathbb{Z})_{torsion}.\] 
%The homology decomposition of $\Sigma Y$ then implies that 
%\[H^{m}(Y_{i};\mathbb{Z})\cong\left\{ 
%     \begin{array}{ll} 
%         H_{m}(\Sigma Y;\mathbb{Z})_{free}\oplus H_{m-1}(\Sigma Y;\mathbb{Z})_{torsion} 
%              & \mbox{if $m\leq i$} \\[2mm]  
%         H_{m}(\Sigma Y;\mathbb{Z})_{torsion} & \mbox{if $m=i+1$} \\[2mm] 
%              0 & \mbox{if $m>i+1$}. 
%      \end{array}\right.\] 
%Thus the map 
%\(\namedright{Y_{i}}{}{\Sigma Y}\) 
%induces an epimorphism in cohomology. 

%Since $\Sigma Y$ is a suspension, all cup products vanish in $H^{\ast}(\Sigma Y;\mathbb{Z})$. 
%Thus, as 
%\(\namedright{H^{\ast}(\Sigma Y;\mathbb{Z})}{}{H^{\ast}(Y_{i};\mathbb{Z})}\) 
%is an epimorphism, the naturality of the cup product implies that all cup products also 
%vanish in $H^{\ast}(Y_{i};\mathbb{Z})$. 
%\end{proof} 

\begin{lemma} 
   \label{W3cup} 
   Cup products vanish in $\widetilde{H}^{\ast}(W_{3};\mathbb{Z}/p^{r}\mathbb{Z})$. 
\end{lemma} 

\begin{proof} 
Recall that $W$ is a $4$-dimensional $CW$-complex with a single $4$-cell.   
Let $Y$ be the 3-skeleton of $W$. Then by cellular approximation and the definition of $W_3$ the inclusion $W_3\hookrightarrow\Sigma W$ factors as a composite
\[
W_3\overset{g}{\to}\Sigma Y\hookrightarrow\Sigma W.
\]
Suppose that there are elements $x,y\in\widetilde{H}^{\ast}(W_{3};\mathbb{Z}/p^{r}\mathbb{Z})$ 
such that $x\cup y\neq 0$. Since $W_{3}$ is simply-connected and of dimension $4$, 
it must be the case that $\vert x\vert=\vert y\vert=2$. By Lemma~\ref{lemma_minimal cell structure}
\[
g^*:H^2(\Sigma Y;\Z/p^r\Z)\to H^2(W_3;\Z/p^r\Z)
\]
is an isomorphism. Let $\bar{x},\bar{y}\in H^2(\Sigma Y;\Z/p^r\Z)$ be elements such that $x=g^{\ast}(\bar{x})$ and $y=g^{\ast}(\bar{y})$. Since $\Sigma Y$ is a suspension, all cup products in $\widetilde{H}^{\ast}(\Sigma Y;\mathbb{Z}/p^{r}\mathbb{Z})$ 
are zero. In particular, we have $\bar{x}\cup\bar{y}=0$. The naturality of the cup product therefore 
implies that
\[
x\cup y=g^{\ast}(\bar{x})\cup g^{\ast}(\bar{y})=g^{\ast}(\bar{x}\cup\bar{y})=0,
\] 
a contradiction. Hence it must be the case that all cup products in 
$\widetilde{H}^{\ast}(W_{3};\mathbb{Z}/p^{r}\mathbb{Z})$ are zero. 
\end{proof} 

\begin{lemma} 
   \label{W3decomp} 
   There is a homotopy equivalence $W_{3}\simeq P\vee\Sigma S\vee\Sigma\overline{P}$.  
\end{lemma}

\begin{proof} 
We will show that the map 
\(\namedright{S\vee\overline{P}}{f_{2}}{W_{2}}\) 
in~(\ref{hlgydecomp}) is null homotopic, implying the statement of the lemma. 
It will be helpful to partition the Moore spaces in $\overline{P}$ by primes. Recall that 
$\overline{P}=\bigvee_{\bar{j}=1}^{\bar{n}} P^{3}(\bar{b}_{\bar{j}})$ where each $\bar{b}_{\bar{j}}$ is an odd prime power. List 
the primes appearing as $\{p_{1},\ldots,p_{t}\}$. Write 
\[\overline{P}=\bigvee_{s=1}^{t}\overline{P}_{s}\qquad\mbox{where}\qquad 
    \overline{P}_{s}=\bigvee_{\ell=1}^{\bar{n}_{s}} P^{3}(p_{s}^{r_{s,\ell}}).\]  
Note that $\bar{n}=\bar{n}_{1}+\cdots +\bar{n}_{t}$. Isolating $\overline{P}_{1}$, let 
\[\overline{Q}=\bigvee_{s=2}^{t} \overline{P}_{s}\] 
so that $\overline{P}=\overline{P}_{1}\vee\overline{Q}$. 
For convenience, write $p_{1}$ as $p$ and $r_{1,\ell}$ as $r_{\ell}$ for $1\leq \ell\leq n_1$ 
so that~$\overline{P}_{1}=\bigvee_{\ell=1}^{n_{1}} P^{3}(p^{r_{\ell}})$. Correspondingly, write 
$P=P_{1}\vee Q$ where $P_{1}$ is the wedge of all the mod-$p^{t}$ Moore spaces in $P$ for 
some $t\geq 1$, and $Q$ is the wedge of mod-$q^{s}$ Moore spaces for all primes $q\neq p$. 
Note that as the torsion in $\overline{P}$ and $P$ may be different, it is possible that for the  
given prime $p$ the wedge $P_{1}$ is trivial. Taking $n_{1}=0$ in the trivial case, write~$P_{1}=\bigvee_{k=1}^{n_{1}} P^{3}(p^{r_{k}})$. The homotopy cofibration 
\(\nameddright{S\vee\overline{P}}{f_{2}}{W_{2}=P}{}{W_{3}}\) 
may then be rewritten as 
\[
S\vee\overline{P}_{1}\vee\overline{Q}\stackrel{f_{2}}{\longrightarrow} P_{1}\vee Q\stackrel{}{\longrightarrow} W_{3}. 
\]

To show that $f_{2}$ is null homotopic it is equivalent to show that each of the composites
\[
\begin{array}{l} 
f_{S}:S\hookrightarrow S\vee\overline{P}_{1}\vee\overline{Q}\overset{f_{2}}{\longrightarrow}P_{1}\vee Q \\[5mm] 
f_{P}:\overline{P}_{1}\hookrightarrow S\vee\overline{P}_{1}\vee\overline{Q}\overset{f_{2}}{\longrightarrow}P_{1}\vee Q\\[5mm]
f_{Q}:\overline{Q}\hookrightarrow S\vee\overline{P}_{1}\vee\overline{Q}\overset{f_{2}}{\longrightarrow}P_{1}\vee Q 
\end{array}
\] 
is null homotopic. Since $f_{2}$ induces the trivial map in integral homology, so do each 
of $f_{P},f_{Q}$ and~$f_S$. 

First, consider $f_{S}$. Since $S$ is $2$-dimensional, $P_{1}\vee Q$ is $1$-connected, 
and $f_{S}$ induces the trivial map in degree two integral homology, the Hurewicz homomorphism implies 
that $f_{S}$ is null homotopic. 

Next, consider $f_{P}$. Since $\overline{P}_{1}=\bigvee_{\ell=1}^{\bar{n}_{1}} P^{3}(p^{r_{\ell}})$, 
to show that $f_{P}$ is null homotopic it suffices to show that the restriction 
\[f^{\ell}_{P}\colon P^{3}(p^{r_{\ell}})\hookrightarrow\namedright{\overline{P}_{1}}{f_{P}}{}{P_{1}\vee Q}\] 
of $f_{P}$ to the $\ell^{th}$ wedge summand is null homotopic. Since $Q$ consists of mod-$q^{s}$ 
Moore spaces for primes $q\neq p$, Lemma~\ref{lemma_[Pa,Pb]=0} implies that $f^{\ell}_{P}$ 
factors as a composite 
\[\namedright{P^{3}(p^{r_{\ell}})}{g_{P}^{\ell}}{P_{1}}\hookrightarrow P_{1}\vee Q\] 
for some map $g_{P}^{\ell}$. We will show that $g_{P}^{\ell}$ is null homotopic, 
thereby implying that $f_{P}^{\ell}$ is null homotopic. 

Observe that as $f_{P}$ induces the zero map in homology, so does $f_{P}^{\ell}$ 
and therefore so does~$g_{P}^{\ell}$. Let $C_{g_{P}^{\ell}}$ be the homotopy cofibre 
of $g_{P}^{\ell}$ and recall that $P_{1}=\bigvee_{k=1}^{n_{1}} P^{3}(p^{r_{k}})$. 
If cup products vanish in $\widetilde{H}^{\ast}(C_{g_{P}}^{\ell};\mathbb{Z}/p^{\min(r_{\ell},r_{k})}\mathbb{Z})$ 
for~$1\leq k\leq n_1$ then Proposition~\ref{trivtestprop} implies that $g_{P}^{\ell}$ is null homotopic. 

It remains to show that cup products vanish in $\widetilde{H}^{\ast}(C_{g_{P}}^{\ell};\mathbb{Z}/p^{\min(r_{\ell},r_k)}\mathbb{Z})$. 
First, as $g_{P}^{\ell}$ induces the zero map in integral homology, by 
Lemma~\ref{lemma_f_* trivial implies all f^* trivial} it also induces the zero map in 
\mbox{mod-$p^{\min(r_{\ell},r_k)}$} cohomology. Second, notice that $g_{P}^{\ell}$ is homotopic to the composite
\[\nameddright{P^{3}(p^{r_{\ell}})}{f_{P}^{\ell}}{P_{1}\vee Q}{pinch}{P_{1}}.\]
The definitions of $f_{P}^{\ell}$ and $f_{P}$ then imply that $g_{P}^{\ell}$ is homotopic to 
the composite 
\[\namedddright{P^{3}(p^{r_{\ell}})}{}{\overline{P}_{1}}{}{S\vee \overline{P}\vee Q}{f_{2}}{P_{1}\vee Q} 
       \stackrel{pinch}{\longrightarrow} P_{1}.\]
As $W_{3}$ is the homotopy cofibre of $f_{2}$ and cup products vanish in 
$\widetilde{H}^{\ast}(W_{3};\mathbb{Z}/p^{\min(r_{\ell},r_k)}\mathbb{Z})$ by Lemma~\ref{W3cup}, 
the factorization of $g_{P}^{\ell}$ through $f_{2}$ and Lemma~\ref{pinchcupabstract} 
imply that cup products vanish in 
$\widetilde{H}^{\ast}(C_{g_{P}}^{\ell};\mathbb{Z}/p^{\min(r_{\ell},r_k)}\mathbb{Z})$.  

Finally, consider $f_{Q}$. Separating out the mod-$p_{s}^{r_{s}}$ Moore spaces 
in $Q$ one prime at a time as was done for $p_{1}$ and $\overline{P}_{1}$, the same argument as 
for $f_{P}$ can be used iteratively. Thus $f_{Q}$ is null homotopic and the proof is complete. 
\end{proof} 

Observe that the space $W_{4}$ in~(\ref{hlgydecomp}) is the same as the suspension 
of the $3$-skeleton of $W$. That is, $W_{4}\simeq\Sigma Y$ for $Y$ the $3$-skeleton 
of $W$. Our approach to dealing with the maps $f_{3}$ and $f_{4}$ in~(\ref{hlgydecomp}) 
will be to use the fact that both $W_{4}$ and $\Sigma W$ are suspensions. This 
requires a general lemma. 

\begin{lemma} 
   \label{retractlemma} 
   Let $A_i$ be simply connected for $1\leq i\leq m$. Suppose that there is a map 
   \(g\colon\namedright{\bigvee_{i=1}^{m} A_{i}}{}{\Sigma X}\) 
   and a sequence $\{i_{1},\ldots,i_{k}\}$ with $1\leq i_{1}<\cdots<i_{k}\leq m$ 
   such that, for~$1\leq j\leq k$, the pinch map 
   \(q_{j}\colon\namedright{\bigvee_{i=1}^{m} A_{i}}{}{A_{i_j}}\) 
   extends across $g$ to a map 
   \(r_{j}\colon\namedright{\Sigma X}{}{A_{i_j}}\). 
   Then the composite 
   \(b\colon\bigvee_{j=1}^{k} A_{i_{j}}\hookrightarrow\namedright{\bigvee_{i=1}^{m} A_{i}}{g}{\Sigma X}\) 
   has a left homotopy inverse. 
\end{lemma} 

\begin{proof} 
Let $r$ be the composite 
\[r\colon\llnameddright{\Sigma X}{\sigma}{\bigvee_{j=1}^{k}\Sigma X}{\bigvee_{j=1}^{k} r_{j}} 
       {\bigvee_{j=1}^{k} A_{i_j}}\] 
where $\sigma$ is defined using the comultiplication on $\Sigma X$. We claim that 
$r\circ b$ is homotopic to a homotopy equivalence. Observe that for $1\leq j\leq k$ we have 
$\tilde{q}_{j}\circ r\simeq r_{j}$ where~\mbox{$\tilde{q}_j:\bigvee_{j=1}^kA_{i_j}\to A_{i_j}$} is the pinch map. By hypothesis, $r_{j}\circ g\simeq q_{j}$, so by definition of $b$ 
we also have $r_{j}\circ b\simeq\tilde{q}_{j}$. Therefore $\tilde{q}_{j}\circ r\circ b\simeq r_{j}\circ b\simeq\tilde{q}_{j}$. 
In homology, the direct sum of finitely many $\mathbb{Z}$-modules is the same as the direct 
product, so the map 
\[\namedright{\widetilde{H}_{\ast}(\bigvee_{j=1}^{k} A_{i};\mathbb{Z})}{r_{\ast}\circ b_{\ast}} 
       {\widetilde{H}_{\ast}(\bigvee_{j=1}^{k} A_{i};\mathbb{Z})\cong\bigoplus_{j=1}^{k} 
       \widetilde{H}_{\ast}(A_{j};\mathbb{Z})}\] 
is determined by the projection to each $\widetilde{H}_{\ast}(A_{j};\mathbb{Z})$. This 
projection is given by $(\tilde{q}_{j})_{\ast}$. Thus the fact that $(\tilde{q}_{j})_{\ast}=(\tilde{q}_{j})_{\ast}\circ r_{\ast}\circ b_{\ast}$ 
implies that $r_{\ast}\circ b_{\ast}$ is the identity map. Hence, by Whitehead's Theorem, 
$r\circ b$ is a homotopy equivalence. 
\end{proof} 

\begin{lemma} 
   \label{W4decomp}
   There is a homotopy equivalence  
   $W_4\simeq P\vee\Sigma S\vee\Sigma\overline{P}\vee\bigvee_{i=1}^{m} S^{4}$.
\end{lemma}

\begin{proof} 
By~(\ref{hlgydecomp}) and Lemma~\ref{W3decomp} there is a homotopy cofibration 
\[
\bigvee_{i=1}^{m} S^3\stackrel{f_{3}}{\longrightarrow} 
    P\vee\Sigma S\vee\Sigma\overline{P}\longrightarrow W_{4} 
\]
where $f_{3}$ induces the trivial map in integral homology. We will show that $f_{3}$ is 
null homotopic and then the statement of the lemma follows. 

Consider the composites 
\begin{equation} 
   \label{f3eqns} 
   \begin{array}{ll} 
       S^{3}\hookrightarrow\namedddright{\bigvee_{i=1}^{m} S^{3}}{f_{3}}{P\vee\Sigma S\vee\Sigma\overline{P}} 
                 {}{P}{}{P^{3}(b_{j})} \\[2mm] 
       S^{3}\hookrightarrow\namedddright{\bigvee_{i=1}^{m} S^{3}}{f_{3}}{P\vee\Sigma S\vee\Sigma\overline{P}} 
                 {}{\Sigma S}{}{S^{3}} \\[2mm] 
       S^{3}\hookrightarrow\namedddright{\bigvee_{i=1}^{m} S^{3}}{f_{3}}{P\vee\Sigma S\vee\Sigma\overline{P}} 
                 {}{\Sigma\overline{P}}{}{P^{4}(\overline{b}_{\bar{j}})}  
    \end{array} 
\end{equation}   
where the three right-hand maps pinch onto a single wedge summand. Let $g$ be the first 
composite in~(\ref{f3eqns}) and let $C_g$ be its cofiber. Since the cofiber of $f_3$ is $W_4$ 
which is the suspension of the 3-skeleton of $W$, all cup products in $\widetilde{H}^{\ast}(W_4;\Z/p_j^{r_j}\Z)$ 
are zero. Therefore, by Lemma~\ref{pinchcupabstract}, all cup products in $\widetilde{H}^{\ast}(C_g;\Z/p_j^{r_j}\Z)$ 
are zero. Hence, by Proposition~\ref{trivtestprop}, $g$ is null homotopic. 

Since $f_{3}$ induces 
the zero map in integral homology, the second and third composites in~(\ref{f3eqns}) are 
null homotopic by the Hurewicz Theorem. These null homotopies hold for the inclusion of each $S^{3}$ 
into $\bigvee_{i=1}^{m} S^{3}$, so $f_{3}$ composes trivially with each of the pinch maps 
\(\namedright{P\vee\Sigma S\vee\Sigma\overline{P}}{}{X}\) 
for $X=P^{3}(b_{j})$, $S^{3}$ or $P^{4}(\bar{b}_{\bar{j}})$. Thus each of these pinch maps extends 
to a map 
\(\namedright{W_{4}}{}{X}\). 
Since $W_{4}$ is a suspension, Lemma~\ref{retractlemma} implies that the map 
\(\namedright{P\vee\Sigma S\vee\Sigma\overline{P}}{}{W_{4}}\) 
has a left homotopy inverse. Hence $f_{3}$ is null homotopic.            
\end{proof} 

\begin{lemma} 
   \label{suspWdecomp} 
   Suppose that $H^{\ast}(W;\mathbb{Z})$ has no $2$-torsion. If the Steenrod operation 
   $Sq^{2}$ acts trivially on $H^{\ast}(W;\mathbb{Z}/2\mathbb{Z})$ then there is a homotopy equivalence 
   \[\Sigma W\simeq P\vee\Sigma S\vee\Sigma\overline{P}\vee\bigg(\bigvee_{i=1}^{m} S^{4}\bigg)\vee S^{5}.\]
   If $Sq^{2}$ acts nontrivially on $H^{\ast}(W;\mathbb{Z}/2\mathbb{Z})$ then there is a 
   homotopy equivalence 
   \[\Sigma W\simeq P\vee\bigvee_{k=2}^{d} S^{3}\vee\Sigma\overline{P}\vee 
           \bigg(\bigvee_{i=1}^{m} S^{4}\bigg)\vee\Sigma\mathbb{C}P^{2}.\] 
\end{lemma} 

\begin{proof} 
By~(\ref{hlgydecomp}) and Lemma~\ref{W4decomp} there is a homotopy cofibration 
\[\nameddright{S^{4}}{f_{4}}{P\vee\Sigma S\vee\Sigma\overline{P}\vee\bigvee_{i=1}^{m} S^{4}}{}{\Sigma W}\]
where $f_{4}$ induces the trivial map in integral homology. Consider the composites 
\begin{equation} 
  \label{f4eqns} 
  \begin{array}{l} 
        \namedddright{S^{4}}{f_{4}}{P\vee\Sigma S\vee\Sigma\overline{P}\vee\bigvee_{i=1}^{m} S^{4}}{}{P} 
              {}{P^{3}(b_{j})} \\[2mm] 
        \namedddright{S^{4}}{f_{4}}{P\vee\Sigma S\vee\Sigma\overline{P}\vee\bigvee_{i=1}^{m} S^{4}}{}{\Sigma S} 
               {}{S^{3}} \\[2mm] 
        \namedddright{S^{4}}{f_{4}}{P\vee\Sigma S\vee\Sigma\overline{P}\vee\bigvee_{i=1}^{m} S^{4}}{} 
              {\Sigma\overline{P}}{}{P^{4}(\bar{b}_{\bar{j}})} \\[2mm] 
        \namedddright{S^{4}}{f_{4}}{P\vee\Sigma S\vee\Sigma\overline{P}\vee\bigvee_{i=1}^{m} S^{4}}{} 
               {\bigvee_{i=1}^{m} S^{4}}{}{S^{4}} 
  \end{array} 
\end{equation} 
where the middle and right maps pinch onto a single wedge summand. 

Suppose that $Sq^{2}$ acts trivially on $H^{\ast}(W;\mathbb{Z}/2\mathbb{Z})$. 
Since each~$b_{j}$ and $\bar{b}_{\bar{j}}$ is a power of an odd prime, by Lemma~\ref{lemma_pi4 P3, pi4 P4 are trivial}, 
$\pi_{4}(P^{3}(b_{j}))\cong\pi_{4}(P^{3}(\bar{b}_{\bar{j}}))\cong 0$ and $\pi_{4}(P^{4}(b_{j}))\cong\pi_{4}(P^{4}(\bar{b}_{\bar{j}}))\cong  0$, implying the first and 
third composites in~(\ref{f4eqns}) are null homotopic. Since 
$\pi_{4}(S^{3})\cong\mathbb{Z}/2\mathbb{Z}$ is generated by a map $\eta$ 
which is detected by $Sq^{2}$, the assumption that $Sq^{2}$ acts trivially on~$H^{\ast}(W;\mathbb{Z}/2\mathbb{Z})$ implies that the second composite 
in~(\ref{f4eqns}) is null homotopic. Since $f_{4}$ induces the zero map in homology, 
the Hurewicz homomorphism implies that the fourth composite in~(\ref{f4eqns}) is null homotopic. 
Thus each of the pinch maps~\mbox{\(\namedright{P\vee\Sigma S\vee\Sigma\overline{P}\vee\bigvee_{i=1}^{m} S^{4}}{}{X}\)}
for $X=P^{3}(b_{j})$, $S^{3}$, $P^{4}(\bar{b}_{\bar{j}})$ or $S^{4}$ extends to a map 
\(\namedright{\Sigma W}{}{X}\). 
Therefore, by Lemma~\ref{retractlemma}, the map 
\(\namedright{P\vee\Sigma S\vee\Sigma\overline{P}\vee\bigvee_{i=1}^{m} S^{4}}{}{\Sigma W}\) 
has a left homotopy inverse. Hence $f_{4}$ is null homotopic, implying that 
\[\Sigma W\simeq P\vee\Sigma S\vee\Sigma\overline{P}\vee\bigg(\bigvee_{i=1}^{m} S^{4}\bigg)\vee S^{5}.\]

Next, suppose that $Sq^{2}$ acts nontrivially on $H^{\ast}(W;\mathbb{Z}/2\mathbb{Z})$. 
Arguing as before, the first, third and fourth composites in~(\ref{f4eqns}) are null 
homotopic. As $Sq^{2}$ detects the generator $\eta$ of $\pi_{4}(S^{3})\cong\mathbb{Z}/2\mathbb{Z}$, 
the nontrivial action of $Sq^{2}$ on $H^{\ast}(W;\mathbb{Z}/2\mathbb{Z})$ implies 
that the second composite in~(\ref{f4eqns}) is nontrivial for at least one of the pinch maps
\(\namedright{\Sigma S=\bigvee_{k=1}^{d} S^{3}}{}{S^{3}\). 
Possibly the second composite in~(\ref{f4eqns}) could be nontrivial for several such 
pinch maps. However, by~\cite{so16}, any map 
\(h\colon\llnamedright{S^{4}}{\bigvee_{k=1}^{d}\epsilon_{k}\eta}{\bigvee_{k=1}^{d} S^{3}}\) 
with $\epsilon_{k}\in\{0,1\}$ for all $1\leq k\leq d$, and having at least one $\epsilon_{k}=1$, 
can be composed with a self-equivalence $e$ of $\bigvee_{k=1}^{d} S^{3}}\) 
so that $e\circ h$ is homotopic to the composite 
\(\namedright{S^{4}}{\eta}{S^{3}}\hookrightarrow\bigvee_{k=1}^{d} S^{3}\) 
where the inclusion can be assumed to be the first wedge summand. Altering 
the copy of $\Sigma S$ in $P\vee\Sigma S\vee\Sigma\overline{P}\vee\bigvee_{i=1}^{m} S^{4}$
by the same self-equivalence $e$, we obtain that each of the pinch maps 
\(\namedright{P\vee\bigvee_{k=2}^{d} S^{3}\vee\Sigma\overline{P}\vee\bigvee_{i=1}^{m} S^{4}}{}{X}\) 
for $X=P^{3}(b_{j})$, $S^{3}$ for $2\leq k\leq d$, $P^{4}(\bar{b}_{\bar{j}})$ or $S^{4}$ extends to a map 
\(\namedright{\Sigma W}{}{X}\). 
Therefore, by Lemma~\ref{retractlemma}, the map 
\(\namedright{P\vee\bigvee_{k=2}^{d} S^{3}\vee\Sigma\overline{P}\vee\bigvee_{i=1}^{m} S^{4}}{}{\Sigma W}\) 
has a left homotopy inverse. Therefore $f_{4}$ factors as the composite 
\[\namedright{S^{4}}{\eta}{S^{3}}\hookrightarrow P\vee\bigvee_{k=1}^{d} S^{3}\vee 
       \Sigma\overline{P}\vee\bigvee_{i=1}^{m} S^{4}\] 
implying that 
$\Sigma W\simeq P\vee\bigvee_{k=2}^{d} S^{3}\vee\Sigma\overline{P}\vee 
           \bigg(\bigvee_{i=1}^{m} S^{4}\bigg)\vee\Sigma\mathbb{C}P^{2}$
since $\Sigma\mathbb{C}P^{2}$ is the homotopy cofibre of $\eta$. 
\end{proof} 

Combining the homotopy decomposition $\Sigma M\simeq\bigg(\bigvee_{i=1}^{m} S^{2}\bigg)\vee\Sigma W$ 
in Lemma~\ref{S1M} with that of $\Sigma W$ in Lemma~\ref{suspWdecomp}, we obtain a homotopy decomposition for $\Sigma M$.

\begin{thm}\label{suspM CW case}
Let $M$ be a 4-dimensional $CW$-complex that has one 4-cell and has homology as in~(\ref{Mhlgy}). If $Sq^2$ acts trivially on $H^*(M;\Z/2\Z)$ then there is a homotopy equivalence 
   \[\Sigma M\simeq\bigg(\bigvee_{i=1}^{\ell} S^{2}\bigg)\vee\bigg(\bigvee_{k=1}^{d} S^{3}\bigg)\vee 
          \bigg(\bigvee_{l=1}^{m} S^{4}\bigg)\vee 
          \bigg(\bigvee_{j=1}^{n} P^{3}(b_{j})\bigg)\vee\bigg(\bigvee_{\bar{j}=1}^{\bar{n}} P^{4}(\bar{b}_{\bar{j}})\bigg)  
          \vee S^{5}.\] 
   If $Sq^2$ acts non-trivially on $H^*(M;\Z/2\Z)$ then there is a homotopy equivalence 
   \[\Sigma M\simeq\bigg(\bigvee_{i=1}^{\ell} S^{2}\bigg)\vee\bigg(\bigvee_{k=1}^{d-1} S^{3}\bigg)\vee 
          \bigg(\bigvee_{l=1}^{m} S^{4}\bigg)\vee 
          \bigg(\bigvee_{j=1}^{n} P^{3}(b_{j})\bigg)\vee\bigg(\bigvee_{\bar{j}=1}^{\bar{n}} P^{4}(\bar{b}_{\bar{j}})\bigg)  
          \vee\Sigma\mathbb{C}P^{2}.\] 
\end{thm} 

As a special case we prove Theorem~\ref{suspM}. 

\begin{proof}[Proof of Theorem~\ref{suspM}]
By assumption $M$ is a smooth, orientable, closed, compact $4$-manifold. Then, by 
Morse Theory, $M$ has a $CW$-structure with one 4-cell. Since $H_{1}(M;\Z)$ is finitely 
generated and has no $2$-torsion, (\ref{Mhyp}) holds and so $H_{\ast}(M;\mathbb{Z})$ 
is as in~(\ref{MPD}). Since~(\ref{MPD}) is a special case of~(\ref{Mhlgy}), Theorem~\ref{suspM CW case} 
applies to decompose $\Sigma M$. Observe that if $M$ is Spin then the Steenrod 
operation $Sq^{2}$ acts trivially on $H^{\ast}(M;\mathbb{Z}/2\mathbb{Z})$, so 
Theorem~\ref{suspM CW case} implies that there is a homotopy equivalence 
\[\Sigma M\simeq\bigg(\bigvee_{i=1}^{m} (S^{2}\vee S^{4})\bigg)\vee 
          \bigg(\bigvee_{j=1}^{n} (P^{3}(b_{j})\vee P^{4}(b_{j}))\bigg)\vee 
          \bigg(\bigvee_{k=1}^{d} S^{3}\bigg)\vee S^{5}, \] 
while if $M$ is non-Spin then $Sq^{2}$ acts nontrivially, so Theorem~\ref{suspM CW case} 
implies that there is a homotopy equivalence 
   \[\Sigma M\simeq\bigg(\bigvee_{i=1}^{m} (S^{2}\vee S^{4})\bigg)\vee 
          \bigg(\bigvee_{j=1}^{n} (P^{3}(b_{j})\vee P^{4}(b_{j}))\bigg)\vee  
          \bigg(\bigvee_{k=1}^{d-1} S^{3}\bigg)\vee\Sigma\mathbb{C}P^{2}.\] 
\end{proof}

\section{Applications} 
\label{sec:apps} 
Suppose that $M$ is a 4-dimensional manifold satisfying the hypotheses of Theorem~\ref{suspM}. In this section 
we give three applications of the homotopy decomposition of $\Sigma M$. 
\medskip 

The first application is to calculate $E^{\ast}(M)$ as a group for any reduced generalized cohomology 
theory $E^{\ast}$. Examples include complex and real $K$-theory and cobordism. 

%A reduced {\color{red} generalized} cohomology theory $h^*$ consists of a sequence of contravariant functors
%\[
%\{h^n:\text{CW}_*\to\text{Ab}\}_{n\in\Z}
%\]
%together with natural isomorphisms $\sigma:h^{n}(X)\overset{\cong}{\to}h^{n+1}(\Sigma X)$ satisfying the following axioms~\cite[Chapter 4E]{hatcher}:
%\begin{itemize}
%\item
%if $f\simeq g:X\to Y$ are based homotopic, then they induce the same homomorphism
%\[
%f^*=g^*:h^n(Y)\to h^n(X);
%\]
%
%\item
%if $A$ is a subcomplex of $X$, then there is an exact sequence
%\[
%h^n(X/A)\to h^n(X)\to h^n(A)
%\]
%
%
%\item
%if $X=\bigvee_{\alpha\in\mathscr{I}}X_{\alpha}$ with inclusions $\imath_{\alpha}:X_{\alpha}\to X$, then there is an isomorphism \[
%\prod_{\alpha\in\mathscr{I}}\imath^*_{\alpha}:h^n(X)\to\prod_{\alpha\in\mathscr{I}}h^n(X_{\alpha}).
%\]
%\end{itemize}

\begin{prop} 
   \label{KM} 
   Let $M$ be a smooth, orientable, closed, connected $4$-manifold satisfying 
   the hypotheses of Theorem~\ref{suspM} and let $E^*$ be a reduced generalized 
   cohomology theory. If $M$ is Spin there is a group isomorphism 
   \[E^n(M)\cong\bigoplus_{i=1}^{m}(E^n(S^{1})\oplus E^n(S^{3}))\oplus\bigoplus_{j=1}^{n}(E^n(P^{2}(b_{j}))\oplus E^n(P^{3}(b_{j}))\oplus\bigoplus_{k=1}^{d}E^n(S^{2})\oplus E^n(S^{4}).\] 
    If $M$ is non-Spin there is a group isomorphism 
   \[E^n(M)\cong\bigoplus_{i=1}^{m}(E^n(S^{1})\oplus E^n(S^{3}))\oplus 
           \bigoplus_{j=1}^{n}(E^n(P^{2}(b_{j}))\oplus E^n(P^{3}(b_{j}))\oplus\bigoplus_{k=2}^{d}E^n(S^{2})\oplus E^n(\mathbb{C}P^{2}).\] 
\end{prop} 

\begin{proof} 
Let $X,A$ and $B$ be CW-complexes such that $\Sigma X\simeq\Sigma A\vee\Sigma B$. Using the axioms of reduced generalized cohomology theories, we obtain a string of group isomorphisms
\begin{eqnarray*}
E^n(X)
&\cong&E^{n+1}(\Sigma X)\\
&\cong&E^{n+1}(\Sigma A\vee\Sigma B)\\
&\cong&E^{n+1}(\Sigma A)\oplus E^{n+1}(\Sigma B)\\
&\cong&E^n(A)\oplus   E^{n}(B)
\end{eqnarray*}
In our case, the asserted group isomorphisms for $E^n(M)$ follow immediately from the above group isomorphisms and the homotopy decomposition of $\Sigma M$ in Theorem~\ref{suspM}.
\end{proof} 

The second application is to current groups. Let $X$ be a smooth manifold and let $G$ be a connected 
Lie group. The \emph{current group} associated to $X$ and $G$ is the space of smooth maps from 
$X$ to~$G$, which is homotopy equivalent to $\map(X,G)$. The most famous example 
is the loop group $\map(S^{1},G)$. Current groups have received considerable attention, 
notably in~\cite{EF,MN,PS}. 

In our case, consider $\map(M,G)$. There is a fibration 
\(\nameddright{\map^{\ast}(M,G)}{}{\map(M,G)}{ev}{G}\) 
where $ev$ evaluates a map at the basepoint of $M$. The multiplication on $G$ 
induces one on~$\map(M,G)$ so the right inverse of $ev$ induced by projecting $M$ 
to the constant map implies that there is a homotopy equivalence 
\begin{equation} 
  \label{mapMGdecomp} 
  \map(M,G)\simeq G\times\map^{\ast}(M,G). 
\end{equation}  
Note that $\map^{\ast}(S^{n},G)=\Omega^{n} G$. For $k\in\mathbb{Z}$, let 
\(\namedright{G}{k}{G}\) 
be the $k^{th}$-power map and let $G\{k\}$ be its homotopy fibre. Applying 
$\map^{\ast}(\ , G)$ to the homotopy cofibration 
\[\nameddright{S^{n}}{k}{S^{n}}{}{P^{n+1}(k)}\] 
gives a homotopy fibration 
\[\nameddright{\map^{\ast}(P^{n+1}(k), G)}{}{\Omega^{n} G}{k}{\Omega^{n} G},\]
implying that $\map^{\ast}(P^{n+1}(k), G)\simeq\Omega^{n} G\{k\}$. 

\begin{prop} 
   \label{current} 
   Let $M$ be a smooth, orientable, closed, connected $4$-manifold satisfying the hypotheses of 
   Theorem~\ref{suspM} and let $G$ be a connected topological group. If $M$ is Spin 
   there is a homotopy equivalence   
   \[\map(M,G)\simeq G\times \prod_{i=1}^{m}(\Omega G\times\Omega^{3} G)\times 
         \prod_{j=1}^{n}(\Omega G\{b_{j}\}\times\Omega^{2} G\{b_{j}\})\times 
         (\prod_{k=1}^{d}\Omega^{2} G)\times\Omega^{4} G.\] 
   If $M$ is non-Spin there is a homotopy equivalence 
   \[\map(M,G)\simeq G\times\prod_{i=1}^{m}(\Omega G\times\Omega^{3} G)\times 
         \prod_{j=1}^{n}(\Omega G\{b_{j}\}\times\Omega^{2} G\{b_{j}\})\times 
         (\prod_{k=2}^{d}\Omega^{2} G)\times\map^{\ast}(\mathbb{C}P^{2},G).\] 
\end{prop} 

\begin{proof} 
In general, if $\Sigma X\simeq\Sigma A\vee\Sigma B$ then 
\[\begin{split} 
    \map^{\ast}(X,G) & \simeq\map^{\ast}(\Sigma X,BG) \\  
        & \simeq\map^{\ast}(\Sigma A,BG)\times\map^{\ast}(\Sigma B,BG) \\  
        & \simeq\map^{\ast}(A,G)\times\map^{\ast}(B,BG). 
  \end{split}\]  
In our case, the homotopy decomposition of $\Sigma M$ in Lemma~\ref{suspM} implies 
that if $M$ is Spin there is a homotopy equivalence  
\[\map^{\ast}(M,G)\simeq\prod_{i=1}^{m}(\Omega G\times\Omega^{3} G)\times 
      \prod_{j=1}^{n}(\Omega G\{b_{j}\}\times\Omega^{2} G\{b_{j}\})\times 
      (\prod_{k=1}^{d}\Omega^{2} G)\times\Omega^{4} G\] 
and if $M$ is non-Spin there is a homotopy equivalence 
\[\map^{\ast}(M,G)\simeq\prod_{i=1}^{m}(\Omega G\times\Omega^{3} G)\times 
      \prod_{j=1}^{n}(\Omega G\{b_{j}\}\times\Omega^{2} G\{b_{j}\})\times 
      (\prod_{k=2}^{d}\Omega^{2} G)\times\map^{\ast}(\mathbb{C}P^{2},G).\]  
The asserted homotopy decompositions for $\map(M,G)$ now follow from~(\ref{mapMGdecomp}). 
\end{proof} 

The third application is to gauge groups. Let $G$ be a simply-connected, simple 
compact Lie group and let $M$ be an orientable, closed, compact $4$-manifold. Then 
$[M,BG]\cong\mathbb{Z}$ so for each~$k\in\mathbb{Z}$ there is a principal 
$G$-bundle $P_{k}$ with second Chern class $k$. The \emph{gauge group}~$\G_{k}(M)$ of $P_{k}$ is the group of $G$-equivariant automorphisms of $P_{k}$ 
that fix $M$. Gauge groups are of paramount importance in mathematical physics 
and geoemetry, and recently their homotopy theory has received a great deal of 
attention~\cite{HK06, HKST18, KK18, KTT17, kono91,KT96, so16, ST18, S, theriault10, 
theriault10b, theriault12, theriault15, theriault17}. 

By~\cite{AB83,gottlieb72} there is a homotopy equivalence 
$B\G_{k}(M)\simeq\map_{k}(M,BG)$ where the right side is the component of 
the space of continuous (not necessarily pointed) maps from $M$ to $BG$ containing 
the map inducing $P_{k}$. From the mapping space point of view there 
is an evaluation fibration sequence 
\[\namedddright{G}{\partial_{k}}{\map^{\ast}_{k}(M,BG)}{}{\map_{k}(M,BG)}{ev}{BG}\] 
where $ev$ evaluates a map at the basepoint of $M$ and $\partial_{k}$ is the fibration 
connecting map. Notice that the homotopy fibre of $\partial_{k}$ is $\G_{k}(M)$. 

In Propositions~\ref{gauge1} and~\ref{gauge2} the Spin and non-Spin cases of smooth, 
orientable, closed, connected $4$-manifolds are considered separately due to some 
additional delicacy in the non-Spin case.

\begin{prop} 
   \label{gauge1} 
   Let $M$ be a smooth, orientable, closed, connected $4$-manifold and let $G$ be a 
   simply-connected, compact, simple Lie group. If $M$ is Spin and satisfies the hypotheses of 
   Theorem~\ref{suspM} then there is a homotopy equivalence 
   \[\G_{k}(M)\simeq\G_{k}(S^{4})\times\prod_{i=1}^{m}(\Omega G\times\Omega^{3} G)\times 
         \prod_{j=1}^{n}(\Omega G\{b_{j}\}\times\Omega^{2} G\{b_{j}\})\times 
         (\prod_{l=1}^{d}\Omega^{2} G).\] 
\end{prop} 

\begin{proof} 
The pinch map 
\(q\colon\namedright{M}{}{S^{4}}\) 
to the top cell induces an isomorphism 
\(\namedright{[S^{4},BG]}{}{[M,BG]}\), 
so by the naturality of the evaluation fibration there is a homotopy fibration diagram 
\begin{equation} 
  \label{natdgrm} 
  \diagram 
       G\rto\ddouble & \map^{\ast}_{k}(S^{4},BG)\rto\dto & \map_{k}(S^{4},BG)\rto^-{ev}\dto & BG\ddouble \\ 
       G\rto & \map^{\ast}_{k}(M,BG)\rto & \map_{k}(M,BG)\rto^-{ev} & BG.  
  \enddiagram 
\end{equation} 
Consider the homotopy cofibration sequence 
\(\namedddright{S^{3}}{f}{M_{3}}{}{M}{q}{S^{4}}\) 
where $M_{3}$ is the $3$-skeleton of $M$ and $f$ is the attaching map for the top cell. This 
induces a homotopy fibration  
\(\nameddright{\map^{\ast}(S^{4},BG)}{}{\map^{\ast}(M,BG)}{}{\map^{\ast}(M_{3},BG)}\). 
Since $\map^{\ast}(M_{3},BG)$ has one component, restricting to the $k^{th}$ component 
of $\map^{\ast}(M,BG)$ we obtain a homotopy fibration 
\(\nameddright{\map^{\ast}_{k}(S^{4},BG)}{}{\map^{\ast}_{k}(M,BG)}{}{\map^{\ast}(M_{3},BG)}\). 
Notice that the connecting map for this homotopy cofibration is $\Sigma f$, which 
is null homotopic by Theorem~\ref{suspM} since it is assumed that $M$ is Spin. 

From the left square in~(\ref{natdgrm}) we therefore obtain a homotopy fibration diagram 
\[\diagram 
      \ast\rto\dto      & \Omega\map^{\ast}_{k}(M,BG)\rdouble\dto & \Omega\map^{\ast}_{k}(M,BG)\dto^{b} \\ 
      \G_{k}(S^{4})\rto\ddouble & \G_{k}(M)\rto^-{a}\dto & \map^{\ast}(\Sigma M_{3},BG)\dto^{(\Sigma f)^{\ast}} \\ 
      \G_{k}(S^{4})\rto\dto & G\rto\dto & \map^{\ast}_{k}(S^{4},BG)\dto \\ 
      \ast\rto	& \map^{\ast}_{k}(M,BG)\rdouble & \map^{\ast}_{k}(M,BG) 
  \enddiagram\] 
where $a$ and $b$ are induced maps. Since $(\Sigma f)^{\ast}$ is null homotopic, 
$b$ has a right homotopy inverse. The homotopy commutativity of the top right square 
then implies that $a$ has a right homotopy inverse. Therefore, using the multiplication 
on $\G_{k}(M)$ we obtain a homotopy equivalence 
\[\G_{k}(M)\simeq\G_{k}(S^{4})\times\map^{\ast}(\Sigma M_{3},BG).\] 
As $M$ is Spin, the homotopy decomposition of $\Sigma M$ in Theorem~\ref{suspM} implies that 
\[\Sigma M_{3}\simeq\bigg(\bigvee_{i=1}^{m} (S^{2}\vee S^{4})\bigg)\vee 
          \bigg(\bigvee_{j=1}^{n} (P^{3}(b_{j})\vee P^{4}(b_{j}))\bigg)\vee 
          \bigg(\bigvee_{l=1}^{d} S^{3}\bigg).\] 
Substituting this into $\map^{\ast}(\Sigma M_{3},BG)$ then gives 
the homotopy equivalence asserted in the statement of the Proposition. 
\end{proof} 

Next, consider the non-Spin case. We aim for an argument mirroring 
the Spin case, but using a map 
\(\namedright{M}{}{\C\PP^{2}}\) 
instead of the pinch map 
\(\namedright{M}{}{S^{4}}\). 
However, the existence of such a map is not obvious. We produce a near substitute 
using the approach in~\cite{so16}. To do so an extra hypothesis is introduced  
on $\pi_{1}(M)$ involving the graph product of groups.

Let $\Gamma=(V, E)$ be a finite undirected graph with vertex 
set $V$ and edge set $E$, and let~\mbox{$\hat{G}=\{G_v|v\in V\}$} be a collection of groups 
associated to the vertices of $\Gamma$. The \emph{graph product $\Gamma\hat{G}$} 
of $\hat{G}$ over $\Gamma$ is the quotient group $F/R$, where $F=*_{v\in V}G_v$ 
is the free product of~$G_v$'s and~$R$ is the normal subgroup generated by 
commutator groups $[G_u,G_v]$ wherever~$(u,v)$ is in $E$. For example, 
if $\Gamma$ is a complete graph then $\Gamma\hat{G}=\bigoplus_{v\in V}G_v$ or 
if $\Gamma$ is a graph of discrete points then $\Gamma\hat{G}=*_{v\in V}G_v$.

If each $G_{v}$ is 
cyclic then the abelianization of $\Gamma\hat{G}$ is $\oplus_{v\in V} G_{v}$. 
It is known that if a group~$H$ is finitely presented then there is a smooth, orientable, closed, 
connected $4$-manifold whose fundamental group is $H$ (see, for 
example, \cite[Theorem 1.2]{F}). For example, if $\Gamma\hat{G}$ is a graph product of cyclic 
groups $\{G_{v}\}_{v\in V}$ then there is a smooth, orientable, closed, connected 
$4$-manifold with 
$\pi_{1}(M)\cong\Gamma\hat{G}$ and $H_{1}(M;\mathbb{Z})\cong\oplus_{v\in V} G_{v}$. 
A specific interesting case is when~$M=M'\times S^{1}$ where $M'$ is a smooth,  
orientable, closed, connected $3$-manifold with~$\pi_{1}(M')$ the graph product 
of copies of $\mathbb{Z}$ (a right-angled Artin group) or copies of $\mathbb{Z}/2\mathbb{Z}$ 
(a right-angled Coxeter group).

\begin{prop} 
   \label{gauge2} 
   Let $M$ be a smooth, orientable, closed, connected $4$-manifold and let $G$ 
   be a simply-connected, compact, simple Lie group. Let $\Gamma\hat{G}$ be a graph  
   product of $\{G_i\}_{i=1}^{m+n}$ where $G_i=\Z$ for $1\leq i\leq m$, $G_{j+m}=\Z/b_j\Z$ 
   for $1\leq j\leq n$, and each $b_j$ is odd. If $M$ is non-Spin and $\pi_1(M)\cong\Gamma\hat{G}$
   then there is a homotopy equivalence 
   \[\G_{k}(M)\simeq\G_{k}(\mathbb{C}P^{2})\times\prod_{i=1}^{m}(\Omega G\times\Omega^{3} G)\times 
         \prod_{j=1}^{n}(\Omega G\{b_{j}\}\times\Omega^{2} G\{b_{j}\})\times 
         (\prod_{l=2}^{d}\Omega^{2} G).\] 
\end{prop} 

\begin{proof} 
For $1\leq i\leq m$, denote the generator of $G_i=\Z$ by $\alpha_i$. For $1\leq j\leq n$, denote the generator of $G_{j+m}=\Z/b_j\Z$ by $\beta_j$. Then each $\alpha_i$ has infinite order and each $\beta_j$ has finite order~$b_j$. Since the Hurewicz homomorphism $h:\pi_1(M)\to H_1(M;\mathbb{Z})$ is the abelianization,~$h(\alpha_i)$ has infinite order and $h(\beta_j)$ has order $b_j$. They generate the direct summands of
\[
H_1(M)\cong\bigoplus^m_{i=1}\Z\oplus\bigoplus^n_{j=1}\Z/b_j\Z.
\]
In particular, $M$ satisfies the hypotheses of Theorem~\ref{suspM}.

For $1\leq i\leq m$, each $\alpha_{i}$ is represented by a map 
\(x_{i}\colon\namedright{S_{1}}{}{M}\) 
of infinite order and for~$1\leq j\leq n$, each $\beta_{j}$ is represented by a map 
\(y_{j}\colon\namedright{S^{1}}{}{M}\) 
of order~$b_{j}$. Since $\beta_j$ has order~$b_j$, it extends to a map 
$\tilde{\beta}_j:P^2(b_j)\to M$. Let 
\[\xi\colon\namedright{\bigg(\bigvee_{i=1}^{m} S^{1}\bigg)\vee\bigg(\bigvee_{j=1}^{n} P^{2}(b_{j})\bigg)} 
      {}{M}\] 
be the wedge sum of the maps $\alpha_{i}$ and $\tilde{\beta}_{j}$. The graph product 
hypothesis on $\pi_{1}(M)$ implies that $\xi$ induces an epimorphism on $\pi_{1}$. By~(\ref{Mhyp}), 
$\xi_{\ast}$ is an isomorphism in degree~$1$ integral homology, and the description 
of $H_{\ast}(M;\mathbb{Z})$ in~(\ref{Mhlgy}) together with the homotopy decomposition 
of $\Sigma M$ in Theorem~\ref{suspM} implies that $\Sigma\xi$ has a left homotopy inverse. 
Define the space $C$ and the map $g$ by the homotopy cofibration 
\[\nameddright{\bigg(\bigvee_{i=1}^{m} S^{1}\bigg)\vee\bigg(\bigvee_{j=1}^{n} P^{2}(b_{j})\bigg)} 
      {\xi}{M}{g}{C}.\] 
Since $\xi$ induces an epimorphism on $\pi_{1}$, $C$ is simply-connected. 
This implies that $C$ can be given a minimal $CW$-structure with one cell 
corresponding to each homology class, and $H_{\ast}(C;\mathbb{Z})$ is determined by 
$H_{\ast}(M;\mathbb{Z})$ since $\zeta_{\ast}$ has a left inverse.  
Since $\Sigma\xi$ has a left homotopy inverse, $\Sigma g$ has a right homotopy inverse. Explicitly, 
the homotopy equivalence for $\Sigma M$ in Theorem~\ref{suspM} implies that 
\[\Sigma C\simeq\bigg(\bigvee_{i=1}^{m} S^{4}\bigg)\vee 
          \bigg(\bigvee_{j=1}^{n} P^{4}(b_{j})\bigg)\vee  
          \bigg(\bigvee_{l=1}^{d-1} S^{3}\bigg)\vee\Sigma\mathbb{C}P^{2}.\] 
This homotopy equivalence may not desuspend but observe that if $C_{3}$ is the $3$-skeleton 
of $C$ then 
\[\Sigma C_{3}\simeq\bigg(\bigvee_{i=1}^{m} S^{4}\bigg)\vee 
          \bigg(\bigvee_{j=1}^{n} P^{4}(b_{j})\bigg)\vee  
          \bigg(\bigvee_{l=1}^{d-1} S^{3}\bigg)\vee S^{3}.\] 
Because $C_{3}$ has cells only in dimensions $2$ and $3$, the attaching maps 
for the $3$-cells are in the stable range, so this homotopy equivalence desuspends and we have
\[
C_3\simeq\bigg(\bigvee_{i=1}^{m} S^{3}\bigg)\vee 
          \bigg(\bigvee_{j=1}^{n} P^{3}(b_{j})\bigg)\vee  
          \bigg(\bigvee_{l=1}^{d-1} S^{2}\bigg)\vee S^2.
\]
Let $D$ be the subwedge of $C_{3}$ given by 
\[D=\bigg(\bigvee_{i=1}^{m} S^{3}\bigg)\vee 
          \bigg(\bigvee_{j=1}^{n} P^{3}(b_{j})\bigg)\vee  
          \bigg(\bigvee_{l=1}^{d-1} S^{2}\bigg).\] 
Then the composite of inclusions 
\(\nameddright{D}{}{C_{3}}{}{C}\) 
has homotopy cofibre $X$, where~\mbox{$\Sigma X\simeq\Sigma\C\PP^{2}$}.

Define the map $q'$ by the composite 
\(q'\colon\nameddright{M}{g}{C}{}{X}\) 
and define the space $Y$ and the maps $f'$ and $\delta$ by the homotopy cofibration sequence
\[
M\overset{q'}{\longrightarrow}X\overset{f'}{\longrightarrow}Y\overset{\delta}{\longrightarrow}\Sigma M 
    \overset{\Sigma q'}{\longrightarrow} \Sigma X. 
\]
As $\Sigma q'$ has a right homotopy inverse 
\(s\colon\namedright{\Sigma X}{}{\Sigma M}\), 
the composite 
\[\nameddright{Y\vee\Sigma X}{\delta\vee s}{\Sigma M\vee\Sigma M}{\nabla}{\Sigma M}\] 
is a homotopy equivalence, where $\nabla$ is the fold map. This implies that $\delta$ has 
a left homotopy inverse and hence $f'$ is null homotopic. Further, when combined with the 
homotopy equivalence for~$\Sigma M$ in Theorem~\ref{suspM}, it implies that there is a 
homotopy equivalence 
\begin{equation} 
  \label{Ydecomp} 
  Y\simeq\bigg(\bigvee_{i=1}^{m} (S^{2}\vee S^{4})\bigg)\vee 
          \bigg(\bigvee_{j=1}^{n} (P^{3}(b_{j})\vee P^{4}(b_{j}))\bigg)\vee 
          \bigg(\bigvee_{l=1}^{d-1} S^{3}\bigg). 
\end{equation}  
Now replace the homotopy cofibration 
\(\nameddright{M}{q}{S^{4}}{\Sigma f}{\Sigma M_{3}}\) 
and the null homotopy for $\Sigma f$ in the argument for the Spin case 
with the homotopy cofibration 
\(\nameddright{M}{}{X}{f'}{Y}\) 
and the null homotopy for $f'$ to obtain a homotopy equivalence 
\[\G_{k}(M)\simeq\G_{k}(X)\times\map^{\ast}(Y,BG).\] 
Substituting the homotopy equivalence for $Y$ in~(\ref{Ydecomp}) into $\map^{\ast}(Y,BG)$ 
then gives a homotopy equivalence  
 \begin{equation} 
   \label{MXsplit} 
   \G_{k}(M)\simeq\G_{k}(X)\times\prod_{i=1}^{m}(\Omega G\times\Omega^{3} G)\times 
         \prod_{j=1}^{n}(\Omega G\{b_{j}\}\times\Omega^{2} G\{b_{j}\})\times 
         (\prod_{l=2}^{d}\Omega^{2} G). 
\end{equation}  
Notice that $X$ only contains one 2-cell and one 4-cell, so it is the cofiber of $a\eta$ for some odd number $a$. While $X$ may not be homotopy equivalent to $\C\PP^{2}$, and while 
$\G_{k}(X)$ may not be homotopy equivalent to $\G_{k}(\C\PP^{2})$, by~\cite[Lemma 2.12]{so16} 
there is a homotopy equivalence~\mbox{$\G_k(X)\times\Omega^2G\simeq\G_k(\C\PP^{2})\times\Omega^2G$} for $d\geq2$. 
If $d=1$, by the construction of $X$, the map $M\to X$ induces isomorphisms 
$H_{\text{free}}^2(M;\Z)\cong H_{\text{free}}^2(X;\Z)$ and $H^4(M:Z)\cong H^4(X;\Z)$. Furthermore, the cup products of degree 2 free elements are preserved under these identifications. So $X$ is a Poincar\'{e} complex and must be $\C\PP^2$. Consequently, $\G_k(X)\simeq \G_k(\C\PP^2)$. Thus, in all cases, from~(\ref{MXsplit})
we obtain the asserted homotopy decomposition of $\G_{k}(M)$.
\end{proof} 
   
Propositions~\ref{gauge1} and~\ref{gauge2} greatly generalize the results in~\cite{so16}, which 
considered the special cases when $\pi_{1}(M)$ is: (i) free, (ii) isomorphic to $\mathbb{Z}/p^{r}\mathbb{Z}$, 
or (iii) a free product of groups in~(i) and~(ii). It is worth emphasizing that the decomposition 
of $\G_k(M)$ can be simply read off from $H_{\ast}(M;\mathbb{Z})$. 

Further, Huang and Wu~\cite{HW17} proved a cancellation result in $p$-local homotopy 
theory. From this we obtain the following. 

\begin{cor}\label{cor_gauge gp of non simply conn reduced to simply conn}
Let $M$ be a manifold as in Propositions~\ref{gauge1} or~\ref{gauge2} and let $p$ be a prime. If $M$ is 
Spin there is a $p$-local homotopy equivalence $\G_k(M)\simeq\G_l(M)$ if and only 
if there is a $p$-local homotopy equivalence $\G_k(S^4)\simeq\G_l(S^4)$. If $M$ 
is non-Spin there is a $p$-local homotopy equivalence $\G_{k}(M)\simeq\G_{l}(M)$ if 
and only if there is a $p$-local homotopy equivalence~\mbox{$\G_k(\C\PP^2)\simeq\G_l(\C\PP^2)$}.~$\qqed$ 
\end{cor}

A classification of when there is a $p$-local homotopy equivalence 
$\G_{k}(S^{4})\simeq\G_{l}(S^{4})$ for any prime $p$ has been determined 
for $G=SU(2)$ \cite{kono91}, $G=SU(3)$ \cite{HK06}, $G=SU(5)$ \cite{theriault15}, 
and~\mbox{$G=Sp(2)$} \cite{theriault10b}. For example, when $G=SU(3)$ there is a 
$p$-local homotopy equivalence $\G_{k}(S^{4})\simeq\G_{l}(S^{4})$ if and only 
if $(k,12)=(l,12)$, where $(a,b)$ is the greatest common denominator of integers 
$a$ and $b$. Partial classifications have been determined in many other cases~\cite{HKST18, 
KK18, KTT17, KT96, ST18, theriault12, theriault17}.

\section*{Acknowledgments}
The first author is funded by PIMS Post-doctoral Fellowship at the University of Regina and was supported by the Oberwolfach Leibniz Fellowship. He would particularly like to thank Mathematisches Forschunginstitut Oberwolfach for the opportunity to work in such a wonderful environment.

\end{document}